\documentclass[english,preprint]{elsarticle}
\usepackage[T1]{fontenc}
\usepackage[latin9]{inputenc}
\usepackage{color}
\usepackage{babel}
\usepackage{prettyref}
\usepackage{amsmath}
\usepackage{amsthm}
\usepackage{amssymb}
\usepackage{graphicx}
\usepackage{xargs}[2008/03/08]
\usepackage[unicode=true,
 bookmarks=true,bookmarksnumbered=true,bookmarksopen=false,
 breaklinks=false,pdfborder={0 0 1},backref=false,colorlinks=true]
 {hyperref}

\makeatletter
  \theoremstyle{definition}
  \newtheorem{defn}{\protect\definitionname}[section]
  \theoremstyle{plain}
  \newtheorem{prop}{\protect\propositionname}[section]
  \theoremstyle{plain}
  \newtheorem{thm}{\protect\theoremname}[section]
  \theoremstyle{plain}
  \newtheorem{cor}{\protect\corollaryname}[section]
 \ifx\proof\undefined\
   \newenvironment{proof}[1][\proofname]{\par
     \normalfont\topsep6\p@\@plus6\p@\relax
     \trivlist
     \itemindent\parindent
     \item[\hskip\labelsep
           \scshape
       #1]\ignorespaces
   }{%
     \endtrivlist\@endpefalse
   }
   \providecommand{\proofname}{Proof}
 \fi
  \theoremstyle{plain}
  \newtheorem{lem}{\protect\lemmaname}[section]

\usepackage{amssymb} %maths
\usepackage{amsmath} %maths
\usepackage{amsfonts}

\usepackage{dimadefs}

\newcommand{\ff}[2]{\ensuremath{(#1)_{#2}}} % Falling factorial

\DeclareMathOperator{\diag}{diag}

\newrefformat{subsec}{Subsection \ref{#1}}
\newrefformat{app}{Appendix \ref{#1}}
\newrefformat{cor}{Corollary \ref{#1}}
\newrefformat{alg}{Algorithm \ref{#1}}
\newrefformat{def}{Definition \ref{#1}}
\newrefformat{clm}{Claim \ref{#1}}
\newrefformat{prop}{Proposition \ref{#1}}
\newrefformat{conj}{Conjecture \ref{#1}}
\newrefformat{prob}{Problem \ref{#1}}
\newrefformat{step}{step \ref{#1}}
\newrefformat{rem}{Remark \ref{#1}}
\newrefformat{tbl}{Table \ref{#1}}
\newrefformat{fact}{Fact \ref{#1}}

\renewcommand{\vec}[1]{\ensuremath{{\boldsymbol #1}}}
\newcommand{\NN}{\ensuremath{\mathbb{N}}}

\newcommand{\CC}{\ensuremath{\mathbb{C}}}

\newcommand{\minsep}{\eta}

\renewcommand{\isdef}{\ensuremath{:=}}

\@ifundefined{showcaptionsetup}{}{%
 \PassOptionsToPackage{caption=false}{subfig}}
\usepackage{subfig}
\makeatother

\providecommand{\corollaryname}{Corollary}
\providecommand{\definitionname}{Definition}
\providecommand{\lemmaname}{Lemma}
\providecommand{\propositionname}{Proposition}
\providecommand{\theoremname}{Theorem}

\begin{document}

\begin{frontmatter}{}

\title{Stability and super-resolution of generalized spike recovery}

\author{Dmitry Batenkov\fnref{fn1,fn2}}
\address{Department of Computer Science, Technion - Israel Institute of Technology, Haifa 32000, Israel}

\fntext[fn1]{Current address: Department of Mathematics, Massachusetts Institute of Technology,
Cambridge, MA 02139, USA. email: batenkov@mit.edu}
\fntext[fn2]{The research leading to these results has received funding from the
European Research Council under European Union's Seventh Framework
Program, ERC Grant agreement no. 320649.}

\ead{batenkov@cs.technion.ac.il}

\begin{abstract}
We consider the problem of recovering a linear combination of Dirac
delta functions and derivatives from a finite number of Fourier samples
corrupted by noise. This is a generalized version of the well-known
spike recovery problem, which is receiving much attention recently.
 We analyze the numerical
conditioning of this problem in two different settings depending on
the order of magnitude of the quantity $N\minsep$, where $N$ is the
number of Fourier samples and $\minsep$ is the minimal distance between
the generalized spikes. In the ``well-conditioned'' regime $N\minsep\gg1$,
we provide upper bounds for first-order perturbation of the solution
to the corresponding least-squares problem. In the near-colliding,
or ``super-resolution'' regime $N\minsep\to0$ with a single cluster, we propose a natural
regularization scheme based on decimating the samples \textendash{}
essentially increasing the separation $\minsep$ \textendash{} and
demonstrate the effectiveness and near-optimality of this scheme in
practice.
\end{abstract}

\end{frontmatter}{}

\global\long\def\nd{\xi}
\global\long\def\jp{z}
\global\long\def\jc{a}
\global\long\def\np{\mathcal{K}}
\newcommandx\meas[1][usedefault, addprefix=\global, 1=k]{m_{#1}}
\global\long\def\df{p}
\global\long\def\ns{N}
\global\long\def\nparams{R}
\global\long\def\fun{f}
\global\long\def\nn#1{\widetilde{#1}}
\global\long\def\err{\varepsilon}

\global\long\def\vec#1{\ensuremath{\mathbf{#1}}}
\global\long\def\cn{\kappa}

\section{Introduction}

In this work we consider the problem of reconstructing the locations
$\nd_{j}\in\left[-\pi,\pi\right]$ and amplitudes $c_{\ell,j}\in\reals$
of a ``generalized spike train''
\begin{align}
f(x) & =\sum_{j=1}^{\np}\sum_{\ell=0}^{\ell_{j}-1}c_{\ell,j}\der{\delta}{\ell}(x-\nd_{j}),\label{eq:gen-delta-fun}
\end{align}
where $\delta\left(x\right)$ is the Dirac delta distribution and
$\der{\delta}{\ell}$ is its derivative of order $\ell$, from a finite
number of the Fourier samples
\[
\widehat{f}\left(k\right)=\frac{1}{2\pi}\int_{-\pi}^{\pi}\fun\left(t\right)\ee^{-\imath kt}\dd t,\quad k=0,1,\dots,\ns-1.
\]

This problem will perhaps be more familiar to the reader in the setting
where $\ell_{j}=1,\;j=1,\dots,\np$, where it becomes the so-called
``spike recovery problem'', receiving much attention recently \cite{azais_spike_2015,candes_super-resolution_2013,candes_towards_2014,donoho_superresolution_1992,demanet_recoverability_2014,dragotti_sampling_2007,fernandez-granda_support_2013,kusuma_accuracy_2009,liao_music_2014,moitra_threshold_2014}.
In this case we have $f\left(x\right)=\sum_{j=1}^{\np}c_{j}\delta\left(x-\nd_{j}\right)$,
and denoting $\jc_{j}:=\frac{c_{j}}{2\pi}$ and $\jp_{j}:=\ee^{-\imath\nd_{j}k}$,
the problem essentially reduces to solving the system of equations
\begin{equation}
\meas=\sum_{j=1}^{\np}\jc_{j}z_{j}^{k},\qquad\jc_{j},z_{j}\in\complexfield,\;k=0,1,\dots,\ns-1.\label{eq:basic-prony}
\end{equation}
This algebraic system appeared originally in the work of G.~R.~de
Prony \cite{prony_essai_1795} in the context of fitting a sum of
exponentials to observed data samples, and hence it is also known as the \emph{Prony system}.
The equations \eqref{eq:basic-prony} appear in areas such as frequency estimation, Padé approximation,
array processing, statistics, interpolation, quadrature, radar signal
detection, error correction codes, and many more \cite{auton_investigation_1981}. 

The higher-order model \eqref{eq:gen-delta-fun} is considered in
many applications, e.g. \cite{batenkov_complete,dragotti_sampling_2007,golub_stable_2000,gustafsson_reconstructing_2000,peter_generalized_2013,sidi_interpolation_1982,van_blaricum_problems_1978}.
In this case, instead of \eqref{eq:basic-prony} we have the following
\emph{polynomial Prony system }(as well as its ``confluent'' variant,
see \cite{batenkov_accuracy_2013}) 
\begin{equation}
\meas=\sum_{j=1}^{\np}z_{j}^{k}\sum_{\ell=0}^{\ell_{j}-1}\jc_{\ell,j}k^{\ell},\qquad\jc_{\ell,j}\in\complexfield,\;\left|z_{j}\right|=1,\label{eq:polynomial-prony}
\end{equation}
where $\jc_{\ell,j}=\frac{c_{\ell,j}\left(-\imath\right)^{\ell}}{2\pi}$.
The unknowns $\left\{ z_{j}\right\} $ (or the corresponding angles
$\nd_{j}=\pm\arg z_{j}$) are frequently called ``poles'', ``nodes''
or ``jumps'', while the linear coefficients $\left\{ \jc_{\ell,j}\right\} $
are called ``magnitudes''.

Issues of numerical stability, or conditioning, of solving \eqref{eq:basic-prony}
and \eqref{eq:polynomial-prony} when the left-hand side is perturbed
have been recognized for a long time. Starting with the original Prony's
method, variety of more stable algorithms have been proposed such
as MUSIC/ESPRIT \cite{rao_model_1992}, matrix pencils \cite{elad_shape_2004,golub_stable_2000},
as well as several least-squares based methods \cite{oleary_variable_2013,pereyra_exponential_2010,peter_nonlinear_2011}
and total variation minimization via convex programming \cite{azais_spike_2015,candes_super-resolution_2013,candes_towards_2014,fernandez-granda_support_2013,moitra_threshold_2014}.
While the majority of these algorithms perform well on simple (i.e. with $\ell_{j}=1$) and
well-separated nodes, they are poorly adapted
to handle either multiple/clustered nodes, non-Gaussian noise or large
values of $\ns$ (\cite{beckermann_numerical_2007,oleary_variable_2013}). An important open problem is \emph{stable super-resolution}, or in other words the possibility to recover closely spaced spikes from noisy measurements, both in \eqref{eq:basic-prony} and all the more in \eqref{eq:polynomial-prony}. Thus in this paper we regard ``super-resolution'' as the regime when the separation is much smaller than $1\over\ns$ \cite{rao_model_1992,stoica_spectral_2005}.

\subsection{Summary of contributions}

In this paper we are mainly interested in the numerical analysis of
the generalized spike recovery problem, and more specifically in understanding the scalings pertaining to the noise amplification. Our first contribution is
providing explicit \emph{component-wise} numerical condition bounds
for the recovery of all the unknown model parameters for the system
\eqref{eq:polynomial-prony}, up to first order, in the overdetermined
setting (i.e. $\ns$ larger than the number of unknowns). Theoretical
analysis of the perturbation for the least-squares solution (\prettyref{sec:Conditioning})
as well as numerical calculations (\prettyref{sec:Numerical-experiments})
of the condition numbers indicate that there is a ``phase transition''
between ill-conditioned and well-conditioned regimes, approximately when the node separation is of the order of $\frac{1}{\ns}$.
Our results describe, in particular, an absolute resolution limit
for any method whatsoever. They build upon and significantly extend
our earlier work \cite{batenkov_accuracy_2013}. 

Our second contribution is proposing a regularization mechanism for
the (mildly) overdetermined Prony problem \eqref{eq:polynomial-prony}
with closely spaced nodes by ``decimation'', i.e. taking subsets of
the equations with indices belonging to arithmetic progressions, and
subsequently solving the resulting square systems (\prettyref{sec:Decimation}).
We show that solution of a decimated system is as accurate as the
(least-squares) solution to the full overdetermined problem (\prettyref{sec:Numerical-experiments}).
Thus, decimation provides a mechanism for achieving near-optimal super-resolution,
at least in the case of a single cluster. 

In \prettyref{sec:Relation-to-existing}, \prettyref{thm:specialization-to-simple-prony}
we specialize the above results to the system \eqref{eq:basic-prony},
and subsequently discuss their relation to existing works in the literature,
in particular \cite{badeau_cramerrao_2008,ben-haim_performance_2012,candes_super-resolution_2013,candes_towards_2014,donoho_superresolution_1992,filbir_problem_2012,kusuma_accuracy_2009,liao_music_2014,maravic2005sar,moitra_threshold_2014,peter_nonlinear_2011,potts_parameter_2010,sidi2003practical}.

\section{\label{sec:Conditioning}Accuracy of the least-squares solution}

\newcommandx\fwm[1][usedefault, addprefix=\global, 1=\ns]{{\cal P}_{#1}}
\global\long\def\CC{\mathbb{C}}
\newcommandx\jac[1][usedefault, addprefix=\global, 1=\ns]{{\cal J}_{#1}}
\newcommandx\pinv[1][usedefault, addprefix=\global, 1=\ns]{\jac[#1]^{\dagger}}
\global\long\def\NN{\mathbb{N}}

\global\long\def\str{\boldsymbol{\ell}}

\subsection{Problem setup}

For a vector $\vec v\in\CC^{m}$, we denote by $\vec v_{k}\;(k=1,\dots,m)$
the $k$-th component of $\vec v$, and we also set $\left|\vec v\right|_{k}:=\left|\vec v_{k}\right|$.
For a matrix $\vec M\in\CC^{m\times n},$ we denote its $i,j$-th
entry by $\vec M_{i,j}$.

In what follows, we assume that the \emph{problem structure vector
$\str=\left(\ell_{1},\dots,\ell_{\np}\right)$ }is fixed\emph{. }We
denote by $\nparams=\nparams\left(\str\right)\isdef\sum_{i=1}^{\np}\ell_{i}+\np$
the overall number of unknown parameters of the problem.

For any $\ns\geqslant\nparams$, we consider the ``forward mapping''
$\fwm:\CC^{\nparams}\to\CC^{\ns}$ given by the measurements \eqref{eq:polynomial-prony},
i.e:
\begin{eqnarray}
\begin{split}\fwm\left(\left(\jc_{0,1},\dots,\jc_{\ell_{1}-1,1},\jp_{1},\dots,\jc_{0,\np},\dots,\jc_{\ell_{\np}-1,\np},\jp_{\np}\right)^{T}\right) & :=\left(\meas[0],\dots,\meas[\ns-1]\right)^{T},\\
\meas & :=\sum_{j=1}^{\np}z_{j}^{k}\sum_{\ell=0}^{\ell_{j}-1}\jc_{\ell,j}k^{\ell}.
\end{split}
\label{eq:fwm-def}
\end{eqnarray}

Thus, we enumerate the $\nparams$ parameters in the order shown \textendash{}
so that $\jc_{0,1}$ is assigned the position $1$, $\jp_{1}$ is
assigned the position $\ell_{1}+1$, and so on. For convenience, we
define
\[
L_{j}:=1+\sum_{m=1}^{j-1}\left(\ell_{m}+1\right),
\]
so that the index corresponding to $\jc_{\ell,j}$ (resp. $z_{j}$)
would be $L_{j}+\ell$ (resp. $L_{j}+\ell_{j}$).

Let $\vec x \in \CC^{\nparams}$ denote a ``data point'' in the parameter space:
\[
\vec x:=\left(\jc_{0,1},\dots,\jc_{\ell_{1}-1,1},\jp_{1},\dots,\jc_{0,\np},\dots,\jc_{\ell_{\np}-1,\np},\jp_{\np}\right)^{T},
\]
so that $\fwm\left(\vec x\right)$ stands for the noise-free measurement vector. The perturbed data vector is $\vec y:=\fwm\left(\vec x\right)+\vec e$.  

Let $\jac\left(\vec x\right)\isdef d\fwm\left(\vec x\right)\in\CC^{\ns\times\nparams}$
denote the Jacobian matrix of the mapping $\fwm$ at the point $\vec x$,
and let $\pinv$ denote the Moore-Penrose pseudo-inverse of $\jac.$

Now consider the solution to the linearized least-squares problem
\begin{equation}\label{eq:linear-lsq-probl}
\vec{x^{*}}=\vec{x^{*}}\left(\vec x, \vec{e}\right):=\arg\min_{\vec z\in\CC^{\nparams}}\left\Vert \vec y-{\cal L}_{\vec x}\left(\vec z\right)\right\Vert ,
\end{equation}
where ${\cal L}_{\vec x}\left(\vec z\right)=\fwm\left(\vec x\right)+\jac\left(\vec x\right)\left(\vec z-\vec x\right)$
is a point in $\left(d\fwm\right)_{\vec x}$, the tangent space of
$\fwm$ at $\vec x$. For small $\|\vec{e}\|$, the vector $\vec{x^*}$ in  \eqref{eq:linear-lsq-probl} is a reasonable proxy for the solution of the nonlinear least squares problem $\vec{x^*}_{nl}:=\arg\min_{\vec{z}\in\CC^{\nparams}} \left\Vert \vec y-\fwm\left(\vec{z}\right)\right\Vert$, and
our main goal in this paper is to investigate the error $\vec{x^{*}}-\vec{x}$. Note that by \eqref{eq:linear-lsq-probl} we have
\[
\vec{x^*} = \arg\min_{\vec{z}\in\CC^{\nparams}} \| \vec{e}-\jac \left( \vec{x} \right) \left( \vec{z}-\vec{x}\right) \|,
\]
and putting $\vec{t}:=\vec{z}-\vec{x}$ this becomes
\begin{equation}\label{eq:ls-pert}
\vec{x^*}-\vec{x}=\arg\min_{\vec{t}\in\CC^{\nparams}} \| \vec{e}-\jac \left( \vec{x} \right)  \vec{t} \|=\pinv\left(\vec x\right)\vec e.
\end{equation}

In order to estimate $\vec{x^{*}}-\vec{x}$, we define the following \emph{component-wise} measure of numerical
conditioning for our problem.
\begin{defn}
\label{def:condition-numbers}Assume that $\jac\left(\vec x\right)$
has full rank. For $\alpha=1,2,\dots,\nparams,$ the \emph{componentwise
condition number }of parameter $\alpha$ at the data point $\vec x\in\CC^{\nparams}$
is the quantity
\begin{equation}
\cn_{\alpha,\ns}\left(\vec x\right):=\sum_{i=1}^{\ns}\left|\pinv\left(\vec x\right)_{\alpha,i}\right|\underbrace{\left|\fwm\left(\vec x\right)_{i}\right|}_{=\left|\meas[i-1]\right|}.\label{eq:cn-def}
\end{equation}
\end{defn}

With this definition, suppose that the measurements have relative error at most $\err$, i.e. that the components of the error vector
$\vec e\in\CC^{\ns}$ satisfy
\begin{equation}
\frac{\left|\vec e_{k}\right|}{\left|\meas\right|}<\err,\quad k=0,1,\dots,\ns-1.\label{eq:noise-model}
\end{equation}
Then by combining \eqref{eq:ls-pert}, \eqref{eq:cn-def} and \eqref{eq:noise-model}, the error of the solution to the linearized least squares problem \eqref{eq:linear-lsq-probl} can be bounded componentwise by
\begin{equation}\label{eq:componentwise-error-bound}
\left|\vec{x^{*}}-\vec x\right|_{\alpha}\leqslant\cn_{\alpha,\ns}\left(\vec x\right)\err.
\end{equation}

In other words, the quantity $\cn_{\alpha,\ns}$ is a measure of noise amplification for the parameter $\alpha$.

The reason for our choice of the noise model \eqref{eq:noise-model} is that the magnitude of the noise-free data \eqref{eq:polynomial-prony} is growing with the index like $m_k \sim k^{\max_j \ell_j-1}$, and so it might be less reasonable to expect the same absolute error in $m_1$ and, say, $m_{100}$ if $\max_j \ell_j > 1$. Other formulations are possible, for instance the absolute error bound $\|\vec{e}_k\|\leqslant \err$, and in fact our results can easily be modified to this scenario\footnote{If instead of the relative noise model \eqref{eq:noise-model} we assume that $\|\vec{e}_k\|\leqslant \err$, we can redefine $\cn_{\alpha,\ns}$ to be just the $\ell_1$ norm of row $\alpha$ of the matrix $\pinv$. The relation \eqref{eq:componentwise-error-bound} would still hold, and the bounds for $\cn_{\alpha,\ns}$ in \prettyref{thm:conditioning-full-stable} would be reduced by the factor $\ns^{\max_j \ell_j-1}$. As a result \emph{all} the parameters of the model \eqref{eq:polynomial-prony} can be stably recovered.  }. However, for reasons of brevity, in the remainder of the paper we shall restrict ourselves to the assumption \eqref{eq:noise-model}.

A central role is played by the \emph{node separation}, defined
as follows.
\begin{defn}
\label{def:delta-zeta}Let $\vec x\in\CC^{\nparams}$ be a data point
such that $\left|z_{j}\right|=1$ for $j=1,\dots,\np.$ For each $j$,
let
\[
\minsep^{\left(j\right)}:=\min_{r\neq j}\left|\arg z_{j}-\arg z_{r}\right|
\]
with the convention that $\minsep^{\left(j\right)}\leqslant\pi$. Furthermore,
we denote 
\[
\minsep=\minsep\left(\vec x\right):=\min_{j}\minsep^{\left(j\right)}.
\]
\end{defn}
\begin{minipage}[t]{1\columnwidth}%
\end{minipage}

Sometimes it will be more convenient to use the absolute distance
instead of the angular distance, i.e.
\begin{equation}
\zeta^{\left(j\right)}:=\min_{r\neq j}\left|\nd_{r}-\nd_{s}\right|,\quad\zeta:=\min_{j}\zeta^{\left(j\right)}\label{eq:zeta-def}
\end{equation}
but clearly 
\begin{equation}
\frac{2}{\pi}\leqslant\frac{\zeta^{\left(j\right)}}{\minsep^{\left(j\right)}},\frac{\zeta}{\minsep}\leqslant1.\label{eq:delta-zeta-corr}
\end{equation}

\begin{minipage}[t]{1\columnwidth}%
\end{minipage}

In what follows, all the constants will in general depend on the problem
structure vector\emph{ $\str$. }Also for consistency we put $\jc_{-1,j}\isdef0$.

Finally, we assume an \emph{a-priori }uniform bound on the magnitudes
of the linear coefficients:
\[
\left|\jc_{\ell,j}\right|\leqslant A.
\]

\subsection{Main result}

It has long been known that the overdetermined Prony system \eqref{eq:basic-prony}
is numerically stable when the number of equations $\ns$ is greater
than $\minsep^{-1}$ (and of course also $\ns\geqslant \nparams$). Here we present a certain quantitative version
of this general principle for the system \eqref{eq:polynomial-prony},
using our definition of condition number as above. For proof see \prettyref{subsec:proof-overdet}.
\begin{thm}
\label{thm:conditioning-full-stable}Let $\vec x\in\CC^{\nparams}$
be a data point such that $\minsep=\minsep\left(\vec x\right)>0$ and
$\jc_{\ell_{j}-1,j}\neq0$ for $j=1,\dots,\np$. Then the Jacobian
matrix $\jac\left(\vec x\right)$
has full rank. Furthermore, there exist constants $K$, $C^{\left(1\right)}$
and $C^{\left(2\right)}$, depending neither on $\ns$ nor $\minsep$,
such that for $\ns>K\cdot\minsep^{-1}$:
\begin{eqnarray*}
\cn_{L_{j}+\ell,\ns}\left(\vec x\right) & \leqslant & C^{\left(1\right)}A\left(1+\frac{\left|\jc_{\ell-1,j}\right|}{\left|\jc_{\ell_{j}-1,j}\right|}\right)\cdot\frac{1}{\ns^{\ell+1-\max_{j}\ell_{j}}},\quad\ell=0,\dots,\ell_{j}-1,\\
\cn_{L_{j}+\ell_{j},\ns}\left(\vec x\right) & \leqslant & C^{\left(2\right)}A\frac{1}{\left|\jc_{\ell_{j}-1,j}\right|}.\frac{1}{\ns^{\ell_{j}+1-\max_{j}\ell_{j}}}.
\end{eqnarray*}
\end{thm}

Note that if the multiplicities of the nodes are different, \prettyref{thm:conditioning-full-stable} shows stability only for the highest-order node. For that node, increasing the number of measurements
$\ns$ improves the accuracy with rate $\frac{1}{\ns}$.
Furthermore, only the highest-order linear coefficient $\jc_{\ell_{j}-1,j}$ is provably
stable, and increasing the number of measurements $\ns$ does not improve
the accuracy for this coefficient beyond a certain bound.
Further note that the (asymptotic) condition numbers themselves \emph{do not depend
on the node separation}, but only the
starting position from which the convergence obeys the stated estimates (the ``well-conditioned'' regime).

In the setting where $\ns\minsep\sim{\cal O}\left(1\right)$, obtaining
comparable estimates for $\cn_{\alpha,\ns}$ appears to be much more
involved. In the remainder of this paper we treat two special cases: the square setting $\ns=\nparams$ (\prettyref{subsec:Square-case}), and a single cluster case (\prettyref{sec:Decimation}).

\subsection{\label{subsec:Square-case}Square case}

For square systems, \prettyref{def:condition-numbers} reduces to
the one used in \cite{batenkov_accuracy_2013}, and in fact it coincides
with the definition of sensitivity of solutions to well-posed algebraic
problems given in \cite{stetter_numerical_2004}. The following estimate
of the conditioning of the system \eqref{eq:polynomial-prony} in
the special case $\ns=\nparams$ is a refinement of the main result
in \cite{batenkov_accuracy_2013}. The proof is presented in \prettyref{subsec:proof-undecimated-refined}.
The main novelty compared to \cite{batenkov_accuracy_2013} is the
explicit dependence on $\minsep^{\left(j\right)}$. 
\begin{thm}
\label{thm:prony-square-cond}Assume that $\ns=\nparams$. Let $\vec x\in\CC^{\nparams}$
be a data point (see \prettyref{def:condition-numbers}) such that
$\minsep=\minsep\left(\vec x\right)>0$ and $\jc_{\ell_{j}-1,j}\neq0$
for $j=1,\dots,\np$. Then the Jacobian matrix $\jac[\nparams]\left(\vec x\right)$
is invertible. Furthermore, there exist constants $C^{\left(3\right)},\;C^{\left(4\right)}$,
not depending on $\minsep$, such that:
\begin{eqnarray*}
\cn_{L_{j}+\ell,\nparams}\left(\vec x\right) & \leqslant & C^{\left(3\right)}A\left(\frac{1}{\minsep^{\left(j\right)}}\right)^{\nparams-\ell}\left(1+\frac{\left|\jc_{\ell-1,j}\right|}{\left|\jc_{\ell_{j}-1,j}\right|}\right),\quad\ell=0,\dots,\ell_{j}-1,\\
\cn_{L_{j}+\ell_{j},\nparams}\left(\vec x\right) & \leqslant & C^{\left(4\right)}A\left(\frac{1}{\minsep^{\left(j\right)}}\right)^{\nparams-\ell_{j}}\cdot\frac{1}{\left|\jc_{\ell_{j}-1,j}\right|}.
\end{eqnarray*}
\end{thm}

\section{\label{sec:Decimation}Decimation}

In contrast with \prettyref{thm:conditioning-full-stable}, now we shift our attention to the ``super-resolution'' regime $\ns \minsep \ll 1$.  In this section we develop a regularization scheme for the 
special case of a single cluster of nodes, based on the idea of \emph{decimation}.
We introduce the \emph{decimated Prony system, }depending on a positive
integer \emph{decimation parameter $\df$, }as follows:
\begin{equation}
n_{k}\isdef\meas[\df k]=\sum_{j=1}^{\np}\jp_{j}^{\df k}\sum_{\ell=0}^{\ell_{j}-1}\left(\jc_{\ell,j}\df^{\ell}\right)k^{\ell},\quad k=0,1,\dots,\nparams-1.\label{eq:decimated-prony-def}
\end{equation}

The idea is that instead of solving \eqref{eq:polynomial-prony} given
$\left\{ \meas[0],\dots,\meas[\ns-1]\right\} $ \textendash{} a difficult
numerical problem \textendash{} one would choose $1\leqslant\df\leqslant\left\lfloor \frac{\ns}{\nparams}\right\rfloor $
and solve the square system \eqref{eq:decimated-prony-def} instead. The main reason why this should work is the following: if we have a cluster of closely spaced nodes $\{\jp_j\}$ with minimal separation $\minsep\ll \pi$, then the modified nodes $\{\jp_j^\df\}$ have minimal separation $\minsep \df$, and therefore by \prettyref{thm:prony-square-cond} these modified nodes can be recovered with improved accuracy by solving \eqref{eq:decimated-prony-def}. Essentially speaking, decimation with parameter $\df$ can be thought of as ``zooming
into'' the cluster by a factor of $\df$. In what follows we provide rigorous justification for this intuition. As it also turns out (see \prettyref{sec:Numerical-experiments}),  the resulting accuracy is
near-optimal, in the sense that it is of the same order as the ``best possible accuracy'' given by the non-decimated
condition number $\cn_{\alpha,\ns}$ in the ``super-resolution'' regime $\ns \minsep \ll 1$. Thus, at least numerically, solving the decimated system provides solution as accurate
as one would get if she solved the full overdetermined problem by
least squares.

\global\long\def\decmap{{\cal P}^{\left(\df\right)}}
\global\long\def\scmap{{\cal R}_{\df}}
\global\long\def\decjac{{\cal J}^{\left(\df\right)}}

Analogously to \prettyref{sec:Conditioning}, we define the decimated
forward map $\decmap:\CC^{\nparams}\to\CC^{\nparams}$ as
\[
\decmap\left(\vec x\right):=\left(n_{0},\dots,n_{\nparams-1}\right),
\]
where $\vec x\in\CC^{\nparams}$ is as in \prettyref{def:condition-numbers}
and $n_{k}$ are given by \eqref{eq:decimated-prony-def}. The decimated
condition numbers $\cn_{\alpha}^{\left(\df\right)}$ are defined as
\[
\cn_{\alpha}^{\left(\df\right)}\left(\vec x\right):=\sum_{i=1}^{\nparams}\left|\left(\left\{ \decjac\left(\vec x\right)\right\} ^{-1}\right)_{\alpha,i}\right|\left|n_{i-1}\right|,
\]
where $\decjac\left(\vec x\right)$ is the Jacobian of the decimated
map $\decmap$ (the definition applies at every point $\vec x$ where
the Jacobian is non-degenerate). In complete analogy to the non-decimated
setting, we set 
\begin{eqnarray*}
\minsep_{\df}^{\left(j\right)}: & = & \min_{r\neq j}\left|\arg z_{r}^{\df}-\arg z_{j}^{\df}\right|,\\
\minsep_{\df} & := & \min_{j}\minsep_{\df}^{\left(j\right)}.
\end{eqnarray*}

The following result is proved in \prettyref{subsec:proof-dec}.
\begin{thm}
\label{thm:decimated-cn-bound}Let $\vec x\in\CC^{\nparams}$ be a
data point (see \prettyref{def:condition-numbers}), and let $\df\geqslant1$
be such that $\minsep_{\df}>0$ and $\jc_{\ell_{j}-1,j}\neq0$ for
$j=1,\dots,\np$. Then the Jacobian matrix $\decjac\left(\vec x\right)$
is invertible. Furthermore, there exist constants $C^{\left(5\right)},\;C^{\left(6\right)}$,
not depending on $\minsep$ and $\df$, such that:
\begin{eqnarray*}
\cn_{L_{j}+\ell}^{\left(\df\right)}\left(\vec x\right) & \leqslant & C^{\left(5\right)}\cdot\left(\frac{1}{\minsep_{\df}^{\left(j\right)}}\right)^{\nparams-\ell}\cdot\left(1+\frac{\left|\jc_{\ell-1,j}\right|}{\left|\jc_{\ell_{j}-1,j}\right|\df^{\ell_{j}-\ell}}\right)\cdot\frac{1}{\df^{\ell+1-\max_{j}\ell_{j}}},\\
\cn_{L_{j}+\ell_{j}}^{\left(\df\right)}\left(\vec x\right) & \leqslant & C^{\left(6\right)}\cdot\left(\frac{1}{\minsep_{\df}^{\left(j\right)}}\right)^{\nparams-\ell_{j}}\cdot\frac{1}{\left|\jc_{\ell_{j}-1,j}\right|}\cdot\frac{1}{\df^{\ell_{j}+1-\max_{j}\ell_{j}}}.
\end{eqnarray*}
\end{thm}
\begin{minipage}[t]{1\columnwidth}%
\end{minipage}

\begin{cor}
\label{cor:super-resolution}Let $\minsep^{*}:=\max_{r\neq j}\left|\arg z_{j}-\arg z_{r}\right|$,
and assume that $\ns\minsep^{*}<\pi\nparams$ (i.e. all nodes form a cluster).
Then the condition numbers of the decimated system \eqref{eq:decimated-prony-def}
with parameter $\df^{*}:=\left\lfloor \frac{\ns}{\nparams}\right\rfloor $
satisfy
\begin{eqnarray*}
\cn_{L_{j}+\ell}^{\left(\df^{*}\right)}\left(\vec x\right) & \leqslant & C^{\left(7\right)}\cdot\left(\frac{1}{\minsep^{\left(j\right)}}\right)^{\nparams-\ell}\left(1+\frac{\left|\jc_{\ell-1,j}\right|}{\left|\jc_{\ell_{j}-1,j}\right|}\right)\cdot\frac{1}{\ns^{\nparams+1-\max_{j}\ell_{j}}},\\
\cn_{L_{j}+\ell_{j}}^{\left(\df^{*}\right)}\left(\vec x\right) & \leqslant & C^{\left(8\right)}\cdot\left(\frac{1}{\minsep^{\left(j\right)}}\right)^{\nparams-\ell_{j}}\frac{1}{\left|\jc_{\ell_{j}-1,j}\right|}\cdot\frac{1}{\ns^{\nparams+1-\max_{j}\ell_{j}}}.
\end{eqnarray*}
\end{cor}
\begin{proof}
Substitution of $\minsep_{\df^{*}}^{\left(j\right)}=\df^{*}\minsep^{\left(j\right)}$
(of course $\minsep_{\df^{*}}^{\left(j\right)}<\pi$) and $\df^{*}:=\left\lfloor \frac{\ns}{\nparams}\right\rfloor $
into \prettyref{thm:decimated-cn-bound} leads to the desired result.
\end{proof}
\begin{minipage}[t]{1\columnwidth}%
\end{minipage}

Comparing \prettyref{cor:super-resolution} with \prettyref{thm:prony-square-cond},
we see an improvement of conditioning by a factor of $\frac{1}{\ns^{\nparams+1-\max_{j}\ell_{j}}}$
(disregarding the constants) gained by decimating - while staying
with the same input size.

Comparing \prettyref{cor:super-resolution} with \prettyref{thm:conditioning-full-stable},
it is seen that if $\minsep$ is fixed, then the decimated condition
numbers (say for the nodes) in the region $\ns\minsep^{*}<\pi\nparams$
decay as $\ns^{-\nparams+\max_{j}\ell_{j}-1}$, while in the region
$\ns\minsep>K$ the rate of decay of the undecimated $\cn$ is only
$\ns^{-\ell_{j}+\max_{j}\ell_{j}-1}$ . This qualitative difference,
or ``phase transition'', is also evident from the numerical data
in \prettyref{sec:Numerical-experiments}.

Let us now discuss how the decimated system can be solved in practice.
\prettyref{cor:super-resolution} provides a simple recipe: given
$\ns$ measurements, just pick up the $\nparams$ evenly spaced ones
having ``maximal spread''. Since this is now a square system (effectively
of constant size), it can be solved efficiently. In \cite{batenkov_prony_2014,batenkov_accurate}
we propose such a method based on polynomial homotopy continuation.
In \prettyref{sec:Numerical-experiments} of this paper we show that
even standard methods such as nonlinear least squares and ESPRIT do
not lose accuracy when provided with decimated measurements on one
hand, and have reduced running time on the other hand.

An important caveat of the decimation approach is that it introduces
\emph{aliasing} for the nodes - indeed, the system \eqref{eq:decimated-prony-def}
has $w_{j}=\jp_{j}^{\df}$ as the solution instead of $\jp_{j}$,
and therefore after solving \eqref{eq:decimated-prony-def}, the algorithm
must select the correct value for the $\df^{\text{-th}}$ root $\left(\nn w_{j}\right)^{\frac{1}{\df}}$.
Thus, either the algorithm should start with an approximation of the
correct value (and thus decimation will be used as a fine-tuning technique),
or it should choose one among the $\df$ possibilities - for instance,
by calculating the discrepancy with the other measurements, which
were not originally utilized in the decimated calculation. Another
possibility would be to try different decimation parameters and employ
some matching procedure, discarding the spurious roots above. In \cite{batenkov_accurate} we discuss these issues in more detail.

\section{Proofs of main results}

\subsection{\label{subsec:proof-undecimated-refined}Proof of \prettyref{thm:prony-square-cond}}

\global\long\def\cpvand{V}
\global\long\def\cvand{U}
\global\long\def\prow#1{\vec{\boldsymbol{v_{#1}}}}
\global\long\def\crow#1{\boldsymbol{u_{#1}}}
\global\long\def\ncoeffs{F}

\global\long\def\powermat#1#2{T_{#1,#2}}
\global\long\def\stirlmat#1{\mathcal{S}_{#1}}

% stirling number of the second kind
%\global\long\def\sts#1#2{\left\{  \begin{array}{c}
% #1\\
% #2 
%\end{array}\right\}  }
\global\long\def\sts#1#2{\frak{S}_{#1}^{(#2)}}

% stirling number of the first kind
%\global\long\def\stf#1#2{\left[\begin{array}{c}
% #1\\
% #2 
%\end{array}\right]}
\global\long\def\stf#1#2{S_{#1}^{(#2)}}

\global\long\def\pbl{E}
\newcommandx\pcoeffbl[1][usedefault, addprefix=\global, 1=j]{\ensuremath{\pbl_{#1}}}

\begin{defn}
\label{def:pvand}Let $\left\{ \ell_{j},\jp_{j}\right\} _{j=1}^{\np}$
be given, and put $\ncoeffs\isdef\sum_{j=1}^{\np}\ell_{j}.$ The\emph{
}Pascal-Vandermonde matrix is the $\ncoeffs\times\ncoeffs$ matrix
\[
\cpvand=\cpvand\left(\jp_{1},\ell_{1},\dots,\jp_{\np},\ell_{\np}\right)\isdef\begin{bmatrix}\prow 0\left(\jp_{1},\ell_{1}\right) & \prow 0\left(\jp_{2},\ell_{2}\right) & \dots & \prow 0\left(\jp_{\np},\ell_{\np}\right)\\
\prow 1\left(\jp_{1},\ell_{1}\right) & \prow 1\left(\jp_{2},\ell_{2}\right) & \dots & \prow 1\left(\jp_{\np},\ell_{\np}\right)\\
\vdots & \vdots & \vdots & \vdots\\
\prow{\ncoeffs-1}\left(\jp_{1},\ell_{1}\right) & \prow{\ncoeffs-1}\left(\jp_{2},\ell_{2}\right) & \dots & \prow{\ncoeffs-1}\left(\jp_{\np},\ell_{\np}\right)
\end{bmatrix},
\]
where 
\[
\prow k\left(\jp_{j},\ell_{j}\right)\isdef\jp_{j}^{k}\left[\begin{array}{ccccc}
1 & k & k^{2} & \dots & k^{\ell_{j}-1}\end{array}\right].
\]
 
\end{defn}
\begin{minipage}[t]{1\columnwidth}%
\end{minipage}
\begin{defn}\label{def:cvand}
Under the above notations, the\emph{ confluent Vandermonde} matrix
is the $\ncoeffs\times\ncoeffs$ matrix
\[
\cvand=\cvand\left(\jp_{1},\ell_{1},\dots,\jp_{\np},\ell_{\np}\right)\isdef\begin{bmatrix}\crow 0\left(\jp_{1},\ell_{1}\right) & \crow 0\left(\jp_{2},\ell_{2}\right) & \dots & \crow 0\left(\jp_{\np},\ell_{\np}\right)\\
\crow 1\left(\jp_{1},\ell_{1}\right) & \crow 1\left(\jp_{2},\ell_{2}\right) & \dots & \crow 1\left(\jp_{\np},\ell_{\np}\right)\\
\vdots & \vdots & \vdots & \vdots\\
\crow{\ncoeffs-1}\left(\jp_{1},\ell_{1}\right) & \crow{\ncoeffs-1}\left(\jp_{2},\ell_{2}\right) & \dots & \crow{\ncoeffs-1}\left(\jp_{\np},\ell_{\np}\right)
\end{bmatrix}
\]
where
\[
\crow k\left(\jp_{j},\ell_{j}\right)\isdef\left[\begin{array}{cccc}
\jp_{j}^{k}, & k\jp_{j}^{k-1}, & \dots & ,\ff{k}{\ell_{j}-1}\jp_{j}^{k-\ell_{j}+1}\end{array}\right]
\]
and $\ff{k}{\ell}$ is the Pochhammer symbol for the falling factorial
\[
\ff{k}{\ell}:=k(k-1)\cdot\dots\cdot(k-\ell+1).
\]
\end{defn}
\begin{minipage}[t]{1\columnwidth}%
\end{minipage}
\begin{defn}
For every $x\in\CC\setminus\left\{ 0\right\} $ and $c\in\naturals$,
let $\powermat xc$ denote the $c\times c$ matrix
\[
\powermat xc:=\diag\left\{ 1,x,x^{2},\dots,x^{c-1}\right\} .
\]
\end{defn}

Clearly,
\[
\left(\powermat xc\right)^{-1}=\powermat{x^{-1}}c.
\]

\begin{minipage}[t]{1\columnwidth}%
\end{minipage}
\begin{defn}
Let $\sts nk$ denote the Stirling number of the second kind \cite[Section 24.1.4]{abramowitz1965handbook}:
\[
\sts nk := \frac{1}{k!} \sum_{j=0}^k \left( -1 \right)^{(k-j)} {k \choose j} j^n,
\]
and let $\stirlmat m$ denote the $m\times m$ upper triangular matrix
\[
\stirlmat m\isdef\left[\begin{array}{ccccc}
\sts 00 & \sts 10 & \sts 20 & \dots & \sts{m-1}0\\
0 & \sts 11 & \sts 21 & \dots & \sts{m-1}1\\
\vdots & \vdots & \ddots &  & \vdots\\
0 & 0 & 0 &  & \sts{m-1}{m-1}
\end{array}\right].
\]
\end{defn}
\begin{minipage}[t]{1\columnwidth}%
\end{minipage}

\begin{prop}
\label{prop:cvand-pvand-conn}The confluent Vandermonde and Pascal-Vandermonde
matrices satisfy
\begin{equation}\label{eq:pascal-vand-confl-conn}
\cpvand\left(\jp_{1},\ell_{1},\dots,\jp_{\np},\ell_{\np}\right)  = \cvand\left(\jp_{1},\ell_{1},\dots,\jp_{\np},\ell_{\np}\right)\times\diag\left\{ \powermat{\jp_{j}}{\ell_{j}}\stirlmat{\ell_{j}}\right\} _{j=1}^{\np}.
\end{equation}
\begin{proof}

The generating function of the Stirling numbers of the second kind is \cite[Section 24.1.4]{abramowitz1965handbook}
\[
\sum_{\ell=0}^{\ell_{j}-1}\sts{\ell_{j}-1}{\ell}\ff{k}{\ell}=k^{\ell}.
\]
The formula \eqref{eq:pascal-vand-confl-conn} then immediately follows from \prettyref{def:pvand} and \prettyref{def:cvand}.
\end{proof}
\end{prop}
\begin{minipage}[t]{1\columnwidth}%
\end{minipage}

The confluent Vandermonde matrix $\cvand$ is well-studied in numerical
analysis due to its central role in polynomial interpolation. The
following fact is well-known \cite{batenkov_accuracy_2013}.
\begin{prop}
\label{prop:cvand-nondegeneracy}The matrix $\cvand\left(\jp_{1},\ell_{1},\dots,\jp_{\np},\ell_{\np}\right)$
is invertible if and only if the nodes $\left\{ \jp_{j}\right\} _{j=1}^{\np}$
are pairwise distinct.
\end{prop}
\begin{minipage}[t]{1\columnwidth}%
\end{minipage}

Now we state the key estimate used to prove \prettyref{thm:prony-square-cond}.
\begin{thm}
\label{thm:vandermonde-inv-bounds}Let $\{x_{1},\dots,x_{n}\}$ be
pairwise distinct complex numbers with $\left|x_{j}\right|\leq1$.
For each $j=1,\dots,n$ assume the separation condition $\left|x_{i}-x_{j}\right|\geq\zeta_{j}>0$
for $i\neq j$. Further, let $\{\ell_{1},\dots,\ell_{n}\}$ be an
ordered collection of natural numbers such that $\ell_{1}+\ell_{2}+\dots+\ell_{n}=N$.
Denote by $\boldsymbol{u}_{j,k}$ the row with index $\ell_{1}+\dots+\ell_{j-1}+k+1$
of $\left[\cvand\left(x_{1},\ell_{1},\dots,x_{n},\ell_{n}\right)\right]^{-1}$
(for $k=0,1,\dots,\ell_{j}-1$). Then the $\ell_{1}$-norm of $\boldsymbol{u}_{j,k}$
satisfies
\begin{equation}
\|\boldsymbol{u}_{j,k}\|_{1}\isdef\sum_{s=1}^{N}\left|\left(\boldsymbol{u}_{j,k}\right)_{s}\right|\leqslant\left(\frac{2}{\zeta_{j}}\right)^{N-\ell_{j}}\frac{2^{k}}{k!}\left(1+\frac{2N}{\zeta_{j}}\right)^{\ell_{j}-1-k}.\label{eq:main-bound}
\end{equation}
\end{thm}
The proof of this theorem (see below) combines original Gautschi's
technique \cite{gautschi1963inverses} and the well-known explicit
expressions for the entries of $\cvand^{-1}$ from \cite{schappelle1972icv},
plus a technical lemma (\prettyref{lem:technical-lemma}).
\begin{defn}
For $j=1,\dots,n$ let
\begin{equation}
h_{j}(x)=\prod_{i\neq j}(x-x_{i})^{-\ell_{i}}.\label{eq:h-def}
\end{equation}
\end{defn}
\begin{lem}
\label{lem:technical-lemma}For any natural $k$, the $k$-th derivative
of $h_{j}$ at $x_{j}$ satisfies
\[
\left|h_{j}^{(k)}\left(x_{j}\right)\right|\leqslant N(N+1)\cdots(N+k-1)\zeta_{j}^{-N-k+\ell_{j}}.
\]
\end{lem}
\begin{proof}
We proceed by induction on $k$. For $k=0$ we have immediately $\left|h_{j}\left(x_{j}\right)\right|\leqslant\zeta_{j}^{-N+\ell_{j}}$.
Now
\begin{equation}
h'_{j}(x)=h_{j}(x)\sum_{i\ne j}\frac{-\ell_{i}}{x-x_{i}}.\label{eq:log-derivative}
\end{equation}
 By the Leibnitz rule we have
\begin{eqnarray*}
h_{j}^{\left(k\right)}\left(x\right) & = & \left(\frac{h_{j}'}{h_{j}}h_{j}\right)^{\left(k-1\right)}\\
 & = & \sum_{r=0}^{k-1}{k-1 \choose r}h_{j}^{(r)}(x)\left(\frac{h'_{j}}{h_{j}}\right)^{\left(k-1-r\right)}\\
 & = & \sum_{r=0}^{k-1}{k-1 \choose r}h_{j}^{(r)}(x)\sum_{i\ne j}\frac{(-1)^{k-r-1}(k-r-1)!\ell_{i}}{(x-x_{i})^{k-r}},
\end{eqnarray*}
hence
\[
\left|h_{j}^{(k)}(x_{j})\right|\leqslant\sum_{r=0}^{k-1}{k-1 \choose r}\left|h_{j}^{(r)}(x_{j})\right|\sum_{i\ne j}\frac{(k-r-1)!\ell_{i}}{|x_{j}-x_{i}|^{k-r}}.
\]
This implies, together with the induction hypothesis, that
\begin{eqnarray*}
\left|h_{j}^{(k)}(x_{j})\right| & \leqslant & \sum_{r=0}^{k-1}{k-1 \choose r}\frac{N\left(N+1\right)\cdots\left(N+r-1\right)}{\zeta_{j}^{N+r-\ell_{j}}}\cdot\frac{(k-r-1)!N}{\zeta_{j}^{k-r}}\\
 & = & \frac{N}{\zeta_{j}^{N+k-\ell_{j}}}\sum_{r=0}^{k-1}\frac{\left(k-1\right)!}{r!}N\left(N+1\right)\cdots\left(N+r-1\right)\\
 & = & \frac{\left(k-1\right)!N}{\zeta_{j}^{N+k-\ell_{j}}}\sum_{r=0}^{k-1}{N-1+r \choose r}.
\end{eqnarray*}
By a well-known binomial identity (proof is immediate by induction
and Pascal's rule) we have
\[
\sum_{r=0}^{k-1}{N-1+r \choose r}={N+k-1 \choose k-1}.
\]
Therefore
\[
\left|h_{j}^{\left(k\right)}\left(x_{j}\right)\right|\leqslant\frac{N\left(N+1\right)\cdots\left(N+k-1\right)}{\zeta_{j}^{N+k-\ell_{j}}},
\]
as required.
\end{proof}
\begin{minipage}[t]{1\columnwidth}%
\end{minipage}
\begin{proof}[Proof of \prettyref{thm:vandermonde-inv-bounds}]
By using a generalization of the Hermite interpolation formula (\cite{spitzbart1960generalization}),
it is shown in \cite{schappelle1972icv} that the components of the
row $\boldsymbol{u}_{j,k}$ are just the coefficients of the polynomial
\[
\frac{1}{k!}\sum_{t=0}^{\ell_{j}-1-k}\frac{1}{t!}h_{j}^{(t)}(x_{j})(x-x_{j})^{k+t}\prod_{i\neq j}(x-x_{i})^{\ell_{i}},
\]
where $h_{j}\left(x\right)$ is given by \eqref{eq:h-def}. By \cite[Lemma]{gautschi1962iva},
the sum of absolute values of the coefficients of the polynomials
$(x-x_{j})^{k+t}\prod_{i\neq j}(x-x_{i})^{\ell_{i}}$ is at most
\[
(1+|x_{j}|)^{k+t}\prod_{i\neq j}(1+|x_{i}|)^{\ell_{i}}\leqslant2^{N-(\ell_{j}-k-t)}.
\]
Therefore
\begin{eqnarray*}
\|\boldsymbol{u}_{j,k}\|_{1} & \leqslant & \frac{1}{k!}\sum_{t=0}^{\ell_{j}-1-k}\frac{1}{t!}\frac{N(N+1)\cdots(N+t-1)}{{\zeta_{j}^{N+t-\ell_{j}}}}2^{N-\ell_{j}+k+t}\\
 & = & \biggl(\frac{2}{\zeta_{j}}\biggr)^{N-\ell_{j}}\frac{2^{k}}{{k!}}\sum_{t=0}^{\ell_{j}-1-k}{\ell_{j}-1-k \choose t}\frac{N(N+1)\cdots(N+t-1)}{(\ell_{j}-k-t)\cdots(\ell_{j}-k-2)(\ell_{j}-k-1)}\biggl(\frac{2}{\zeta_{j}}\biggr)^{t}\\
 & \leqslant & \biggl(\frac{2}{\zeta_{j}}\biggr)^{N-\ell_{j}}\frac{2^{k}}{{k!}}\biggl(1+\frac{2N}{\zeta_{j}}\biggr)^{\ell_{j}-1-k},
\end{eqnarray*}
which completes the proof (in the last transition we used $\frac{N+r}{s+r}\leqslant\frac{Ns+rN}{s+r}=N$,
where $s=\ell_{j}-k-t\geqslant1$ and $r=0,\dots,t-1$).
\end{proof}
\begin{minipage}[t]{1\columnwidth}%
\end{minipage}

Now we state a similar bound for the Pascal-Vandermonde matrix $\cpvand$.
\begin{cor}
\label{cor:rowwise-norm-cvand-2}Assume that $\left|\jp_{j}\right|=1$,
with $\min_{i\neq j}\left|\jp_{i}-\jp_{j}\right|=\zeta_{j}>0$ for
$j=1,\dots,\np$. Denote by $\boldsymbol{v}_{j,k}$ the row with index
\begin{equation}
\ell_{1}+1+\dots+\ell_{j-1}+1+k+1\label{eq:index-in-cor-norm-cvand-2}
\end{equation}
 of $\left\{ \cpvand\left(\jp_{1},\ell_{1}+1,\dots,\jp_{\np},\ell_{\np}+1\right)\right\} ^{-1}$
(for $k=0,1,\dots,\ell_{j}$). Then there exists a constant $C$,
not depending on $\zeta$, such that 
\begin{equation}
\|\boldsymbol{v}_{j,k}\|_{1}\leqslant C\cdot\left(\frac{1}{\zeta_{j}}\right)^{\nparams-k}\label{eq:cvand-norm-bound-2}
\end{equation}
where $\nparams=\sum_{j=1}^{\np}\left(\ell_{j}+1\right)=\ncoeffs+\np$.
\end{cor}
\begin{proof}
Denote by $u_{j,k}$ the row with index \eqref{eq:index-in-cor-norm-cvand-2}
of
\[
\left\{ \cvand\left(\jp_{1},\ell_{1}+1,\dots,\jp_{\np},\ell_{\np}+1\right)\right\} ^{-1}.
\]
 Since $\left(\powermat{\jp_{j}}{\ell_{j}}\stirlmat{\ell_{j}}\right)^{-1}$
is block upper triangular with entries bounded by a constant\footnote{As a matter of fact, we have the exact formula for the inverse \cite[Section 24.1.4]{abramowitz1965handbook}
\[
\stirlmat m^{-1}=\left[\begin{array}{ccccc}
\stf 00 & \stf 10 & \stf 20 & \dots & \stf{m-1}0\\
0 & \stf 11 & \stf 21 & \dots & \stf{m-1}1\\
\vdots & \vdots & \ddots &  & \vdots\\
0 & 0 & 0 &  & \stf{m-1}{m-1}
\end{array}\right],
\]
where $\stf nk$ is the Stirling number of the first kind, equal to the (signed) number of permutations of $n$ sybmols having exactly $k$ cycles \cite[Section 24.1.3]{abramowitz1965handbook}  .}, say, $C^{*}$, we have by \prettyref{thm:vandermonde-inv-bounds}
(obviously $\zeta_{j}<2\nparams)$
\begin{eqnarray*}
\|\boldsymbol{v}_{j,k}\|_{1} & \leqslant & \ell_{j}\cdot C^{*}\cdot\max_{t=k,\dots,\ell_{j}}\|u_{j,t}\|_{1}\\
 & \leqslant & C^{*}\ell_{j}\left(\frac{2}{\zeta_{j}}\right)^{\nparams-\ell_{j}}\max_{t=k,\dots,\ell_{j}}\frac{2^{t}}{t!}\left(\frac{4\nparams}{\zeta_{j}}\right)^{\left(\ell_{j}+1\right)-1-t}\\
 & \leqslant & C\cdot\left(\frac{1}{\zeta_{j}}\right)^{\nparams-k},
\end{eqnarray*}
which finishes the proof.
\end{proof}
\begin{minipage}[t]{1\columnwidth}%
\end{minipage}
\begin{defn}
For every $j=1,\dotsc,\np$ let us denote by $\pcoeffbl$ the following
$(\ell_{j}+1)\times(\ell_{j}+1)$ block
\begin{eqnarray}
\pcoeffbl=\pcoeffbl\left(\vec x\right) & \isdef & \left[\begin{array}{cccc}
1 & 0 & 0 & 0\\
0 & 1 & 0 & \frac{\jc_{0,j}}{\jp_{j}}\\
\vdots & \vdots & \ddots & \vdots\\
0 & 0 & 0 & \frac{\jc_{\ell_{j}-1,j}}{\jp_{j}}
\end{array}\right].\label{eq:coeffbl-def}
\end{eqnarray}

Subsequently, we denote by $\pbl$ the block diagonal $\nparams\times\nparams$
matrix
\begin{equation}
\pbl=\pbl\left(\vec x\right):=\diag\left\{ \pcoeffbl[1],\dots,\pcoeffbl[\np]\right\} .\label{eq:pbl-def}
\end{equation}
\end{defn}
\begin{prop}
Direct calculation gives
\begin{align}
\pcoeffbl^{-1} & =\begin{bmatrix}1 & 0 & 0 & \dots & 0\\
0 & 1 & 0 & \dots & -\frac{\jc_{0,j}}{\jc_{\ell_{j}-1,j}}\\
0 & 0 & 1 & \dots & -\frac{\jc_{1,j}}{\jc_{\ell_{j}-1,j}}\\
 &  &  &  & \vdots\\
0 & 0 & 0 & \dots & +\frac{\jp_{j}}{\jc_{\ell_{j}-1,j}}
\end{bmatrix},\label{eq:pcoeffbl-inv}
\end{align}
\end{prop}
\begin{proof}[Proof of \prettyref{thm:prony-square-cond}]

For the Jacobian matrix of $\fwm[\nparams]$, we have the following
factorization by a straightforward computation:
\begin{equation}
\jac[\nparams]\left(\vec x\right)=\cpvand\left(\jp_{1},\ell_{1}+1,\dots,\jp_{\np},\ell_{\np}+1\right)\times\pbl\left(\vec x\right).\label{eq:nondec-jac-fact}
\end{equation}
Therefore 
\begin{eqnarray*}
\jac[\nparams]^{-1} & = & \diag\left\{ \pcoeffbl[j]^{-1}\right\} \cpvand^{-1}.
\end{eqnarray*}
Combining this with \eqref{eq:pcoeffbl-inv}, \eqref{eq:delta-zeta-corr}
and \prettyref{cor:rowwise-norm-cvand-2}, we complete the proof of
\prettyref{thm:prony-square-cond}.
\end{proof}

\subsection{\label{subsec:proof-dec}Proof of \prettyref{thm:decimated-cn-bound}}

From \eqref{eq:decimated-prony-def} it is clear that the map $\decmap$
can be written as a composition: $\decmap=\fwm[\nparams]\circ\scmap$,
where $\fwm[\nparams]$ is given by \eqref{eq:fwm-def} and $\scmap$
is the rescaling mapping given by
\begin{eqnarray*}
\scmap\left(\left(\jc_{0,1},\dots,\jc_{\ell_{1}-1,1},\jp_{1},\dots,\jc_{0,\np},\dots,\jc_{\ell_{\np}-1,\np},\jp_{\np}\right)^{T}\right) & \isdef\\
\left(b_{0,1},\dots,b_{\ell_{1}-1,1},w_{1},\dots,b_{0,\np},\dots,b_{\ell_{\np}-1,\np},w_{\np}\right)^{T} & =\\
\left(\jc_{0,1}\cdot\df^{0},\dots,\jc_{\ell_{1}-1,1}\cdot\df^{\ell_{1}-1},\jp_{1}^{\df},\dots,\jc_{0,\np}\cdot\df^{0},\dots,\jc_{\ell_{\np}-1,\np}\cdot\df^{\ell_{\np}-1},\jp_{\np}^{\df}\right)^{T}.
\end{eqnarray*}
By the chain rule, $d\decmap=d\fwm[\nparams]\times d\scmap$. But
$d\scmap$ is just the diagonal matrix
\[
d\scmap=\diag\left\{ 1,\df,\df^{2},\dots,\df^{\ell_{1}-1},\df\jp_{1}^{\df-1},\dots,1,\df,\df^{2},\dots,\df^{\ell_{\np}-1},\df\jp_{\np}^{\df-1}\right\} .
\]
By definition, $\min_{r\neq j}\left|\arg w_{r}-\arg w_{j}\right|=\minsep_{\df}^{\left(j\right)}$.
Furthermore, we have the estimate $\left|n_{k}\right|\leqslant A\df^{\max_{j}\ell_{j}-1}$.
Taking the inverse, and applying \prettyref{cor:rowwise-norm-cvand-2}
and \eqref{eq:nondec-jac-fact}, it can be seen that the decimated
condition numbers satisfy:
\begin{eqnarray*}
\cn_{L_{j}+\ell}^{\left(\df\right)}\left(\vec x\right) & \leqslant & C^{\left(5\right)}A\left(\frac{1}{\minsep_{\df}^{\left(j\right)}}\right)^{\nparams-\ell}\left(1+\frac{\left|b_{\ell-1,j}\right|}{\left|b_{\ell_{j}-1,j}\right|}\right)\cdot\frac{1}{\df^{\ell+1-\max_{j}\ell_{j}}},\\
\cn_{L_{j}+\ell_{j}}^{\left(\df\right)}\left(\vec x\right) & \leqslant & C^{\left(6\right)}\left(\frac{1}{\minsep_{\df}^{\left(j\right)}}\right)^{\nparams-\ell_{j}}\frac{1}{\left|b_{\ell_{j}-1}\right|}\cdot\frac{1}{\df^{2-\max_{j}\ell_{j}}}.
\end{eqnarray*}
Now plug in $\left|b_{\ell,j}\right|=\df^{\ell}\left|\jc_{\ell,j}\right|$
to finish the proof of \prettyref{thm:decimated-cn-bound}.

\subsection{\label{subsec:proof-overdet}Proof of \prettyref{thm:conditioning-full-stable}}

\global\long\def\cvcol{\vec w}

\global\long\def\pvandf{W}
\global\long\def\ww{{\cal W}}

A key step in the proof of this result is an accurate description
of pseudo-inverses of rectangular Pascal-Vandermonde\emph{ }matrices,
with the nodes on the unit circle. 
\begin{defn}
Let $\left\{ \ell_{j},\jp_{j}\right\} _{j=1}^{\np}$ be given. For
any $t=0,1,\dots,$ and $j=1,\dots,\np,$ denote by $\cvcol_{j,\ns}^{\left(t\right)}$
the column vector (where $0^{0}=1$) 
\[
\cvcol_{j,\ns}^{\left(t\right)}=\begin{pmatrix}0^{t}\\
z_{j}\\
2^{t}z_{j}^{2}\\
\vdots\\
\left(\ns-1\right)^{t}z_{j}^{\ns-1}
\end{pmatrix}.
\]
With this notation, we define the following $\ns\times\nparams$ matrix:
\[
\pvandf_{\ns}=\pvandf_{\ns}\left(\jp_{1},\ell_{1},\dots,\jp_{\np},\ell_{\np}\right):=\begin{pmatrix}\cvcol_{1,\ns}^{\left(0\right)} & \dots & \cvcol_{1,\ns}^{\left(\ell_{1}\right)} & \dots & \cvcol_{\np,\ns}^{\left(0\right)} & \dots & \cvcol_{\np,\ns}^{\left(\ell_{\np}\right)}\end{pmatrix}.
\]

We also put
\[
\ww:=\pvandf_{\ns}^{*}\pvandf_{\ns}\in\CC^{\nparams\times\nparams}.
\]
\end{defn}
Recalling \prettyref{def:pvand}, note that $\pvandf_{\nparams}=\cpvand\left(\jp_{1},\ell_{1}+1,\dots,\jp_{\np},\ell_{\np}+1\right)$.
Thus we immediately obtain the following corollary of \prettyref{prop:cvand-nondegeneracy}.

\begin{minipage}[t]{1\columnwidth}%
\end{minipage} 
\begin{prop}
\label{prop:full-rank}Suppose that $\left\{ \jp_{j}\right\} $ are
pairwise distinct. Then, for $\ns\geqslant\nparams,$ the matrix $\pvandf_{\ns}$
has full column rank, and thus $\ww$ has full rank.
\end{prop}
\begin{minipage}[t]{1\columnwidth}%
\end{minipage} 

The next claim is easily verified by observation.
\begin{prop}
\label{prop:vhv-block-structure}The matrix $\ww$ has an explicit
block structure as follows:
\[
\ww=\left[B_{rs}\right]_{1\leqslant r,s\leqslant\np},
\]
where $B_{rs}$ is a rectangular $\left(\ell_{r}+1\right)\times\left(\ell_{s}+1\right)$
block
\[
B_{rs}=\left[b_{i,j}^{\left(r,s\right)}\right]_{0\leqslant i\leqslant\ell_{r},\;0\leqslant j\leqslant\ell_{s}}
\]
and 
\begin{equation}
b_{i,j}^{\left(r,s\right)}=\left[\cvcol_{r,\ns}^{\left(i\right)}\right]^{*}\cvcol_{s,\ns}^{\left(j\right)}=\sum_{\ell=0}^{\ns-1}\ell^{i+j}\left(z_{r}^{*}z_{s}\right)^{\ell}.\label{eq:b-entries}
\end{equation}
\end{prop}
\begin{minipage}[t]{1\columnwidth}%
\end{minipage}
\begin{defn}
Given $N,q$ integers, $h_{N,q}$ is the sum of $q$-th \emph{powers
(generalized harmonic number)
\[
h_{N,q}\isdef\sum_{\ell=0}^{N-1}\ell^{q}.
\]
}
\end{defn}
For instance, $h_{\ns,0}=\ns$, $h_{\ns,1}=1+\dots+\left(\ns-1\right)=\frac{\ns\left(\ns-1\right)}{2}$.
In general, by the Faulhaber's formula \cite{conway_book_1995}, $h_{\ns,q}$ is a polynomial
in $\ns$ with leading term $\frac{1}{q+1}\ns^{q+1}$.
\begin{prop}
\label{prop:b-ijrs-entries-mag}The entries $b_{i.j}^{\left(r,s\right)}$
defined in \eqref{eq:b-entries} satisfy, as $\ns>K_{1}$ for some
constant $K_{1}$ (depending only on $i+j$) 
\[
\left|b_{i,j}^{\left(r,s\right)}\right|\leqslant\begin{cases}
\frac{2}{i+j+1}\ns^{i+j+1} & r=s,\\
4\minsep_{rs}^{-1}\ns^{i+j} & r\neq s,
\end{cases}
\]
where $\minsep_{rs}:=\left|\arg z_{r}^{*}z_{s}\right|$ (as in \prettyref{def:delta-zeta}).
\end{prop}
\begin{proof}
Let $u_{rs}=z_{r}^{*}z_{s}$. It is a complex number on the unit circle. Consider two cases.

\begin{enumerate}
\item $r=s$ and so $u_{rs}=1$. In this case $b_{i,j}^{\left(r,r\right)}$
is just the $\left(i+j\right)$-th generalized harmonic number 
\[
b_{i,j}^{\left(r,r\right)}=h_{\ns,i+j}.
\]
\item $r\neq s$. Let $q$ be a non-negative integer, put  $u_{rs}:=z$ and consider
\[
f_{\ns,q}\left(z\right):=\sum_{k=0}^{\ns-1}k^{q}z^{k}.
\]
We evaluate the above expression using summation by parts. Define
sequences $A_{k}:=k^{q}$ and 
\[
B_{k}:=1+z+\dots+z^{k-1}.
\]
That is, $B_{k+1}-B_{k}=z^{k}$ with $B_{0}:=0$. Thus
\begin{eqnarray*}
f_{\ns,q}\left(z\right) & = & \sum_{k=0}^{\ns-1}A_{k}\left(B_{k+1}-B_{k}\right)\\
 & = & A_{\ns}B_{\ns}-A_{0}B_{0}-\sum_{k=0}^{\ns-1}B_{k+1}\left(A_{k+1}-A_{k}\right)\\
 & = & \ns^{q}B_{\ns}-\sum_{k=0}^{\ns-1}\left[\left(k+1\right)^{q}-k^{q}\right]B_{k}.
\end{eqnarray*}
Now put $z=\exp\left(\imath t\right)$ (without loss of generality
for $0<t<\pi$). Then obviously for any non-negative integer $k$
we have 
\begin{eqnarray*}
\left|B_{k}\right|^{2} & = & \left|\frac{z^{k+1}-1}{z-1}\right|^{2}=\left|\frac{\sin\left(k+1\right)\frac{t}{2}}{\sin\frac{t}{2}}\right|^{2},
\end{eqnarray*}
and thus $\left|B_{k}\right|\leqslant\frac{2}{t}$. Therefore
\begin{eqnarray*}
\left|f_{\ns,q}\left(z\right)\right| & \leqslant & \frac{2}{t}\left\{ \ns^{q}+\sum_{k=0}^{\ns-1}\left[\left(k+1\right)^{q}-k^{q}\right]\right\} =\frac{4}{t}\ns^{q}.
\end{eqnarray*}
\end{enumerate}
This proves the claim.
\end{proof}
\begin{minipage}[t]{1\columnwidth}%
\end{minipage}

Now we move on to study $\ww^{-1}$.
\begin{prop}
\label{prop:b-rr-inv-entries}The square matrix $B_{rr}$ is invertible,
with $\left(i,j\right)$-th entry ($i,j$ starting from $1$) of the
inverse satisfying for $\ns>K_{2}$
\[
\left(B_{rr}^{-1}\right)_{i,j}\leqslant C_{1}\cdot\frac{q_{i,j}}{\ns^{i+j-1}},
\]
where $q_{i.j}$ is the $\left(i.j\right)$-th entry of the inverse
$\left(\ell_{r}+1\right)\times\left(\ell_{r}+1\right)$ Hilbert matrix,
and $C_{1}$, as well as $K_{2}$, do not depend on $\ns$.
\end{prop}
\begin{proof}
Use formula for component-wise perturbation of matrix inverse. Namely,
write
\[
B_{rr}=H_{\ell_{r}}+\Delta H
\]
where $H_{\ell_{r}}$ is the scaled $\left(\ell_{r}+1\right)\times\left(\ell_{r}+1\right)$
Hilbert matrix
\begin{equation}
H_{\ell_{r}}=\left(\frac{\ns^{i+j-1}}{i+j-1}\right)_{i,j}.\label{eq:h-mag}
\end{equation}

Given any matrix $A$, let us denote by $\left|A\right|$ the matrix
of absolute values of entries of $A$. Now we have $\left|\Delta H\right|\leqslant\epsilon\left|H_{\ell_{r}}\right|$
for $\epsilon\sim\ns^{-1}$. It is immediately checked that
\begin{equation}
H_{\ell_{r}}^{-1}=\left(\frac{q_{i,j}}{\ns^{i+j-1}}\right)_{i,j}\label{eq:h-inv-mag}
\end{equation}
where $q_{i,j}$ is the $\left(i.j\right)$-th entry of the inverse
$\left(\ell_{r}+1\right)\times\left(\ell_{r}+1\right)$ Hilbert matrix.

Then (see \cite[Section 3]{higham_survey_1994}) to first order in
$\epsilon$ we have $B_{rr}^{-1}=H_{\ell_{r}}^{-1}+\Delta B_{rr}^{-1}$
where
\[
\left|\Delta B_{rr}^{-1}\right|\sim\left|H_{\ell_{r}}^{-1}\right|\left|H_{\ell_{r}}\right|\left|H_{\ell_{r}}^{-1}\right|\epsilon.
\]
Taking into account the order of magnitudes specified by \eqref{eq:h-mag}
and \eqref{eq:h-inv-mag} we easily obtain that the order of growth
of $\left(B_{rr}^{-1}\right)_{i,j}$ is 
\[
\frac{q_{i,j}}{\ns^{i+j-1}}+O\left(\ns^{-i-j}\right).
\]
Since the entries of $B_{rr}$ are polynomials in $\ns$ (see \prettyref{prop:b-ijrs-entries-mag}),
the entries of $B_{rr}^{-1}$ are rational functions in $\ns$, and
thus we obtain the desired result.
\end{proof}
\begin{minipage}[t]{1\columnwidth}%
\end{minipage}

Now we come to the main structure result for $\ww$. 

\begin{minipage}[t]{1\columnwidth}%
\end{minipage}
\begin{defn}
Given the structure vector $\str=\left(\ell_{1},\dots,\ell_{\np}\right)$,
let $D_{\str}$ denote the following block diagonal matrix:
\[
D_{\str}=\diag\left\{ B_{11},\dots,B_{\np\np}\right\} .
\]
\end{defn}
Recall that the matrix $\ww$ consists of the rectangular blocks $B_{rs}$.
The following claim is straightforward.
\begin{prop}
We have
\[
\ww=D_{\vec{\ell}}\times X,
\]
where $X\in\CC^{\nparams\times\nparams}$ has the block structure
\[
X=\left[C_{rs}\right]_{1\leqslant r,s\leqslant\np},
\]
each $C_{rs}$ being a $\left(\ell_{r}+1\right)\times\left(\ell_{s}+1\right)$
block
\[
C_{rs}=B_{rr}^{-1}\times B_{rs}.
\]
So in particular $C_{rr}=I_{\left(\ell_{r}+1\right)\times\left(\ell_{r}+1\right)}.$ 
\end{prop}
\begin{minipage}[t]{1\columnwidth}%
\end{minipage}

Now using \prettyref{prop:b-ijrs-entries-mag} and \prettyref{prop:b-rr-inv-entries}
we easily obtain the following.
\begin{prop}
For $r\neq s$, the $\left(i,j\right)$-th entry of $C_{rs}$ (counting
starts from $1$) satisfies, for $\ns>K_{2}$ and some constant $C_{2}$
\[
\left|\left[C_{rs}\right]_{i,j}\right|\leqslant C_{2}\cdot\minsep^{-1}\ns^{-i+j-1}.
\]
\end{prop}
\begin{minipage}[t]{1\columnwidth}%
\end{minipage}

Next we denote $Y:=I_{\nparams\times\nparams}-X$. By induction on
$k$, it is easy to prove the following fact.
\begin{prop}
For each $k=1,2,\dots,$ the matrix $Y^{k}$ has the block structure
\[
Y^{k}=\left[T_{rs}^{\left(k\right)}\right]_{1\leqslant r,s\leqslant\np},
\]
where $T_{rs}^{\left(k\right)}$ is a $\left(\ell_{r}+1\right)\times\left(\ell_{s}+1\right)$
block, whose $\left(i.j\right)$-th entry satisfies, for $\ns>K_{2}$
and some constant $C_{3}$
\[
\left|\left[T_{rs}^{\left(k\right)}\right]_{i,j}\right|\leqslant C_{3}\cdot\frac{\nparams^{k-1}}{\minsep^{k}}\ns^{-i+j-k}.
\]
\end{prop}
This immediately leads to the following conclusion.
\begin{prop}
For $\ns>K_{3}:=\max\left(\frac{\nparams}{\minsep},K_{2}\right)$ the
Neumann series $\sum_{k=1}^{\infty}Y^{k}$ converges, and thus $X=I-Y$
is invertible, with
\[
X^{-1}=I+\sum_{k=1}^{\infty}Y^{k}=I+Z,
\]
where $Z$ has the same block structure as $X$, i.e. $Z=\left[\Xi_{rs}\right]_{1\leqslant r,s\leqslant\np}$
, with $\Xi_{rs}$ being a $\left(\ell_{r}+1\right)\times\left(\ell_{s}+1\right)$
block, whose $\left(i,j\right)$-th entry satisfies, for some constant
$C_{4}$ 
\[
\left|\left[\Xi_{rs}\right]_{i,j}\right|\leqslant C_{4}\cdot\frac{1}{1-\frac{\nparams}{\ns\minsep}}\cdot\begin{cases}
\ns^{-i+j-1} & r\neq s,\\
\ns^{-i+j-2} & r=s
\end{cases}.
\]
\end{prop}
Now, since $\ww=D_{\vec{\ell}}\left(I-Y\right)$, then
\begin{eqnarray*}
\ww^{-1} & = & X^{-1}D_{\vec{\ell}}^{-1}=\left(I+Z\right)D_{\vec{\ell}}^{-1}\\
 & = & D_{\vec{\ell}}^{-1}+\left[\Xi_{rs}\right]\diag\left\{ B_{tt}^{-1}\right\} .
\end{eqnarray*}
Using all the above structural results, we obtain the following asymptotic
description of the blocks of $\ww^{-1}$.
\begin{prop}
The matrix $\ww^{-1}\in\CC^{\nparams\times\nparams}$ has the block
form 
\[
\ww^{-1}=\left[{\cal V}_{rs}\right]_{1\leqslant r,s\leqslant\np},
\]
where each ${\cal V}_{rs}$ is a $\left(\ell_{r}+1\right)\times\left(\ell_{s}+1\right)$
block, whose $\left(i,j\right)$-th entry satisfies, for some constant
$C_{5}$ and $\ns>K_{3}$,
\[
\left|\left[{\cal V}_{rs}\right]_{i,j}\right|\leqslant C_{5}\cdot\frac{1}{1-\frac{\nparams}{\ns\minsep}}\cdot\begin{cases}
\ns^{-i-j+1} & r=s,\\
\ns^{-i-j} & r\neq s.
\end{cases}
\]
\end{prop}
So we actually have proved the following result.
\begin{thm}
\label{thm:vand-pinv-magn}Consider the pseudo-inverse $\pvandf_{\ns}^{\dagger}=\ww^{-1}\pvandf_{\ns}^{*}\in\CC^{\nparams\times\ns}$
Pascal-Vandermonde matrix as a striped matrix, i.e. $\pvandf_{\ns}^{\dagger}=\left[\vec v_{\ell,j}\right]_{1\leqslant j\leqslant\np}^{0\leqslant\ell\leqslant\ell_{j}}$,
where each $\vec v_{\ell,j}\in\CC^{1\times\ns}$ is a row vector.
Then as $\ns>K_{4}:=\max\left(K_{3},\frac{2\nparams}{\minsep}\right)$,
the magnitudes of the entries of $\vec v_{\ell,j}$ are bounded by
$C_{6}\cdot\ns^{-\ell-1}$, where $C_{6}$ depends only on the problem
structure vector $\str$.
\end{thm}
\begin{minipage}[t]{1\columnwidth}%
\end{minipage}
\begin{proof}[Proof of \prettyref{thm:conditioning-full-stable}]

For the Jacobian matrix $\jac\left(\vec x\right)=d\fwm\left(\vec x\right)\in\CC^{\ns\times\nparams}$,
direct computation gives

\[
\jac\left(\vec x\right)=\pvandf_{\ns}\times\pbl,
\]
where $\pbl$ is defined in \eqref{eq:pbl-def}. Combining this with
\prettyref{prop:full-rank} proves that $\jac$ has
full rank.

Furthermore, 
\begin{eqnarray*}
\pinv=\left(\jac^{*}\jac\right)^{-1}\jac^{*} & = & \left(\pbl^{*}\pvandf_{\ns}^{*}\pvandf_{\ns}\pbl\right)^{-1}\pbl^{*}\pvandf_{\ns}^{*}=\pbl^{-1}\ww^{-1}\left(\pbl^{*}\right)^{-1}\pbl^{*}\pvandf_{\ns}^{*}\\
 & = & \pbl^{-1}\ww^{-1}\pvandf_{\ns}^{*}=\pbl^{-1}\pvandf_{\ns}^{\dagger}.
\end{eqnarray*}

Consider $\jac^{\dagger}\in\CC^{\nparams\times\ns}$ as a striped
matrix, i.e. $\jac^{\dagger}=\left[\vec j_{\ell,j}\right]_{1\leqslant j\leqslant\np}^{0\leqslant\ell\leqslant\ell_{j}}$
where each $\vec j_{\ell,j}\in\CC^{1\times\ns}$ is a row vector.
Using \eqref{eq:pcoeffbl-inv} and \prettyref{thm:vand-pinv-magn},
we obtain that for $\ns>K_{4}$ and some constant $C_{7}$ 
\begin{equation}
\left|\left(\vec j_{\ell,j}\right)_{t}\right|\leqslant C_{7}\cdot\begin{cases}
\left(1+\frac{\left|\jc_{\ell-1,j}\right|}{\left|\jc_{\ell_{j}-1,j}\right|}\right)\cdot\frac{1}{\ns^{\ell+1}} & 0\leqslant\ell<\ell_{j},\\
\frac{1}{\left|a_{\ell_{j}-1,j}\right|}\cdot\frac{1}{\ns^{\ell_{j}+1}} & \ell=\ell_{j}.
\end{cases}\label{eq:jac-pinv-entries}
\end{equation}

Let $i=0,1,\dots,\ns-1$. Clearly, we have
\[
\left|\fwm\left(\vec x\right)\right|_{k}=\left|\meas[k-1]\right|\leqslant C_{8}A\left(k-1\right)^{\max_{j}\ell_{j}-1}.
\]

Thus in particular
\begin{equation}
\sum_{k=0}^{\ns-1}\left|m_{k}\right|\leqslant C_{9}A\ns^{\max_{j}\ell_{j}}.\label{eq:moment-abs-vals-sum}
\end{equation}
Plugging \eqref{eq:jac-pinv-entries} and \eqref{eq:moment-abs-vals-sum}
into \eqref{eq:cn-def}, the second claim of \prettyref{thm:conditioning-full-stable}
immediately follows.
\end{proof}

\section{\label{sec:Numerical-experiments}Numerical experiments}

\subsection{Condition numbers}

In this section we present numerical study of the quantities $\cn_{\alpha,\ns}$
and $\cn_{\alpha}^{\left(\df^{*}\right)}$, and their comparison with
the respective upper bounds given by \prettyref{thm:conditioning-full-stable}
and \prettyref{cor:super-resolution}.

\subsubsection{Experimental setup}
\begin{enumerate}
\item In all experiments, the nodes were chosen to be evenly spaced and
of the same order (i.e. $\ell_{r}=\ell_{s}=n$ for all $r,s$). In
all the experiments we put $\np=3$. The variable parameters were
$n$ and $\minsep$. 
\item We were interested primarily in asymptotics w.r.t $\ns$ and $\minsep$.
Thus, in order to minimize the influence of the magnitudes of the
linear coefficients $\jc_{\ell,j}$, we effectively computed the inner
products of the rows of the corresponding (pseudo-) inverse Vandermonde
matrices $W_{\ns}^{\dagger}$ and $V^{-1}$ with the measurement vector,
see \prettyref{subsec:proof-undecimated-refined} and \prettyref{subsec:proof-overdet}.
\item The following quantities were computed:

\begin{enumerate}
\item Decimated and undecimated condition numbers.
\item Theoretical bounds for the stable regime, according to \prettyref{thm:conditioning-full-stable}
(accurate computation was done according to \prettyref{prop:b-rr-inv-entries},
and specifically \eqref{eq:h-mag} and \eqref{eq:h-inv-mag}):
\[
Bound1_{\ell}\left(N\right):=\ns^{n}\left(H_{\ell_{r}}^{-1}\right)_{\ell+1,1}.
\]
\item Theoretical bounds for the super-resolution regime, according to \prettyref{cor:super-resolution}
(see also \prettyref{cor:rowwise-norm-cvand-2}):
\[
Bound2_{\ell}\left(\ns,\minsep\right):=\frac{\nparams^{\nparams}2^{\left(\nparams+2\left(\ell_{j}-\ell\right)+1\right)}\cdot\ell_{j}}{\ell!}\minsep^{\ell-\nparams}\ns^{n-1-\nparams}.
\]
\end{enumerate}
\item All calculations were done using Mathematica with 30 digit precision.
\end{enumerate}

\subsubsection{Results}

\begin{figure}
\subfloat[$n=3,\;\minsep=0.1,\;\protect\np=3$]{\includegraphics[width=0.9\textwidth]{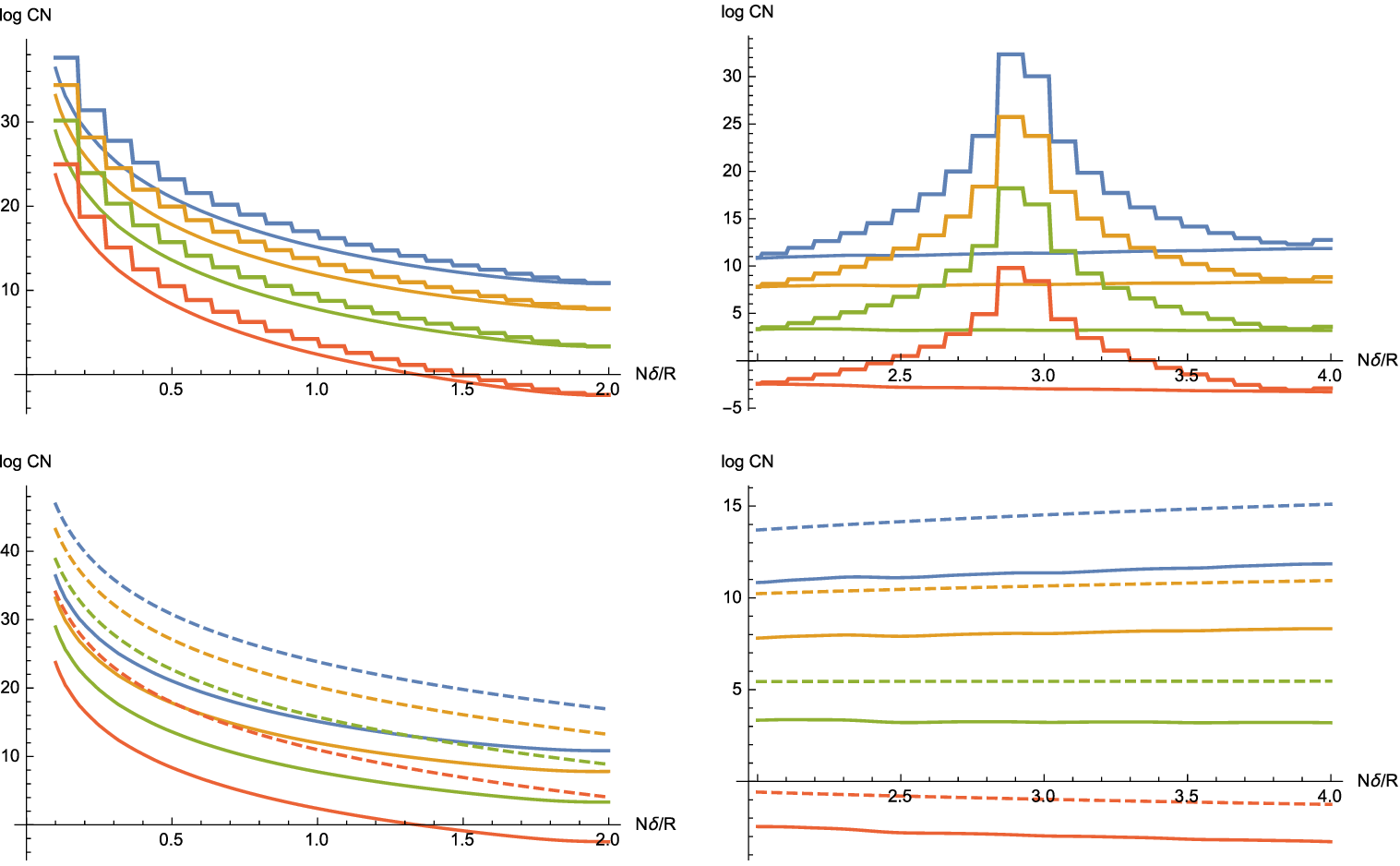}

}

\subfloat[$n=2,\;\minsep=0.05,\;\protect\np=3$]{\includegraphics[width=0.9\textwidth]{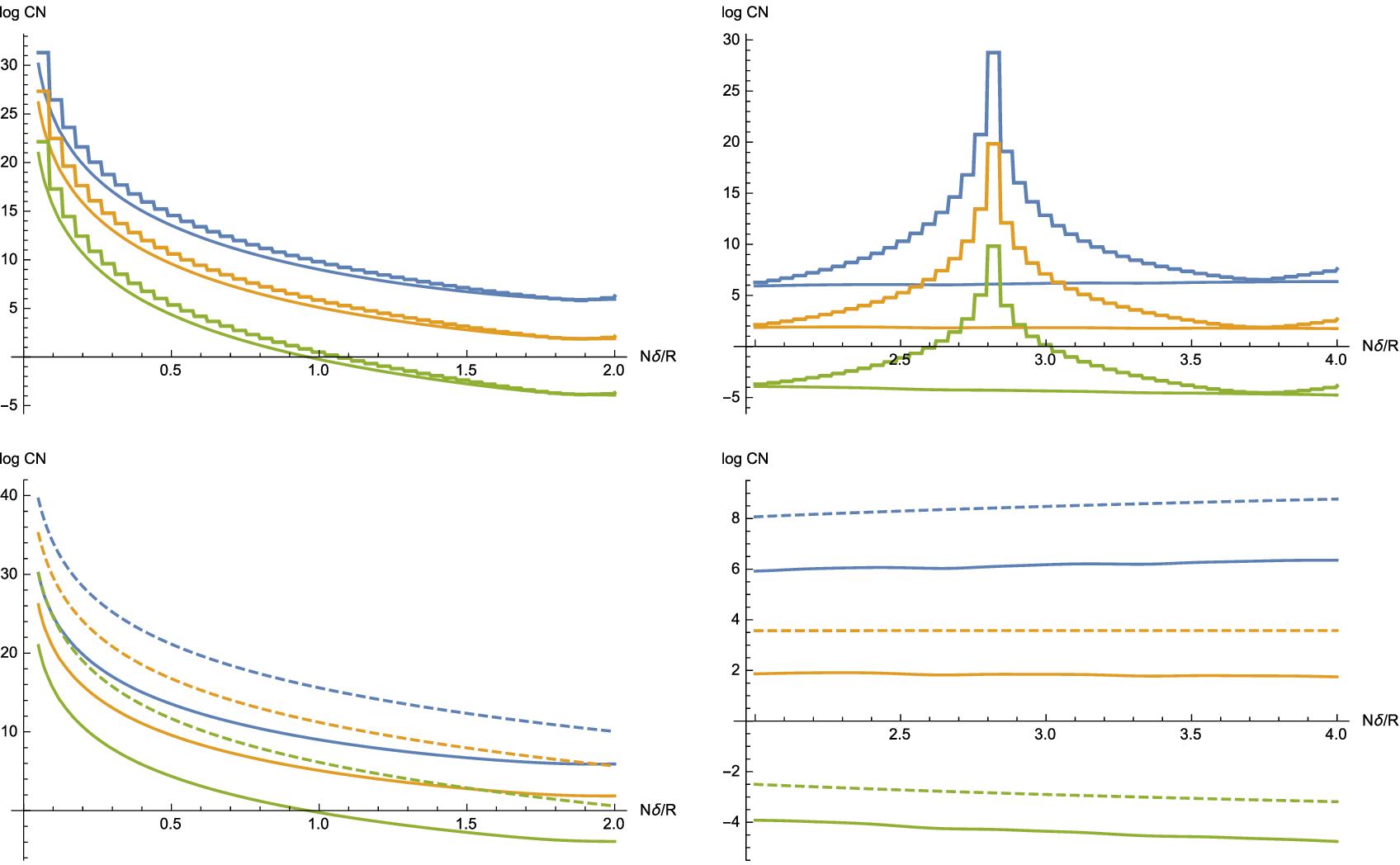}}

\caption{Estimating the condition numbers and their upper bounds. Upper row:
decimated vs. non-decimated, super-resolution (upper left) and stable
(upper right) regions. Lower row: undecimated condition numbers vs.
upper bounds. CN stands for condition number ($\protect\cn$) and $\delta$ stands for $\minsep$.}
\label{fig:cond-numb}

\end{figure}

The graphs in \prettyref{fig:cond-numb} present the computed values
of $\cn_{L_{j}+\ell,\ns}$ (solid) and $\cn_{L_{j}+\ell}^{\left(\df^{*}\right)}$
(thick solid), as well as the quantities $Bound1_{\ell}\left(\ns\right)$
(dashed) and $Bound2_{\ell}\left(\ns,\minsep\right)$ (dotted). The
different values of $\ell$ are distinguished by color-coding. In
each experiment we fixed $\np$, $n$ and $\minsep$, while varying
$\ns$. The horizontal axis is scaled as $\frac{\ns\minsep}{\nparams}$.
The plots are semi-logarithmic in the vertical axis.

\subsubsection{Conclusions}
\begin{enumerate}
\item A ``phase transition'' between well-conditioned and ill-conditioned
regions is seen to occur with the threshold in the range $\frac{\ns\minsep}{\nparams}\in\left(1,3\right)$.
\item In the ``near ill-conditioned'' (or ``super-resolution'') region,
the decimated condition number are almost identical with the non-decimated
ones.
\item The computed upper bounds provide accurate growth rates in the region
$\ns\minsep\gg1$, and are also relatively accurate in the super-resolution
region. 
\item The periodic pattern for $\cn^{\left(\df\right)}$ is seen in the
well-conditioned region and it is well-predicted by the theory. For
instance, it is easy to see that for infinite number of values of
$\df$ we have $\pi<\df\minsep^{*}<\pi+\varepsilon$ (recall \prettyref{cor:super-resolution}),
thus $\minsep_{\df}$ becomes small and $\cn^{\left(\df\right)}$ blows
up.
\end{enumerate}

\subsection{Least Squares and ESPRIT with decimation}

We have tested the decimation technique on two well-known algorithms
for Prony systems - generalized ESPRIT \cite{badeau2008performance}
and nonlinear least squares (LS, implemented by MATLAB's \texttt{lsqnonlin}).
To avoid the aliasing problem, we assumed an initial approximation
to be given. All computations were done in MATLAB with double precision
floating point arithmetic. The computed values of $\meas$ were perturbed
in a random manner with specified noise level. 

\begin{figure}
\subfloat[ESPRIT, $\ell_{j}=3$]{\includegraphics[width=0.45\columnwidth]{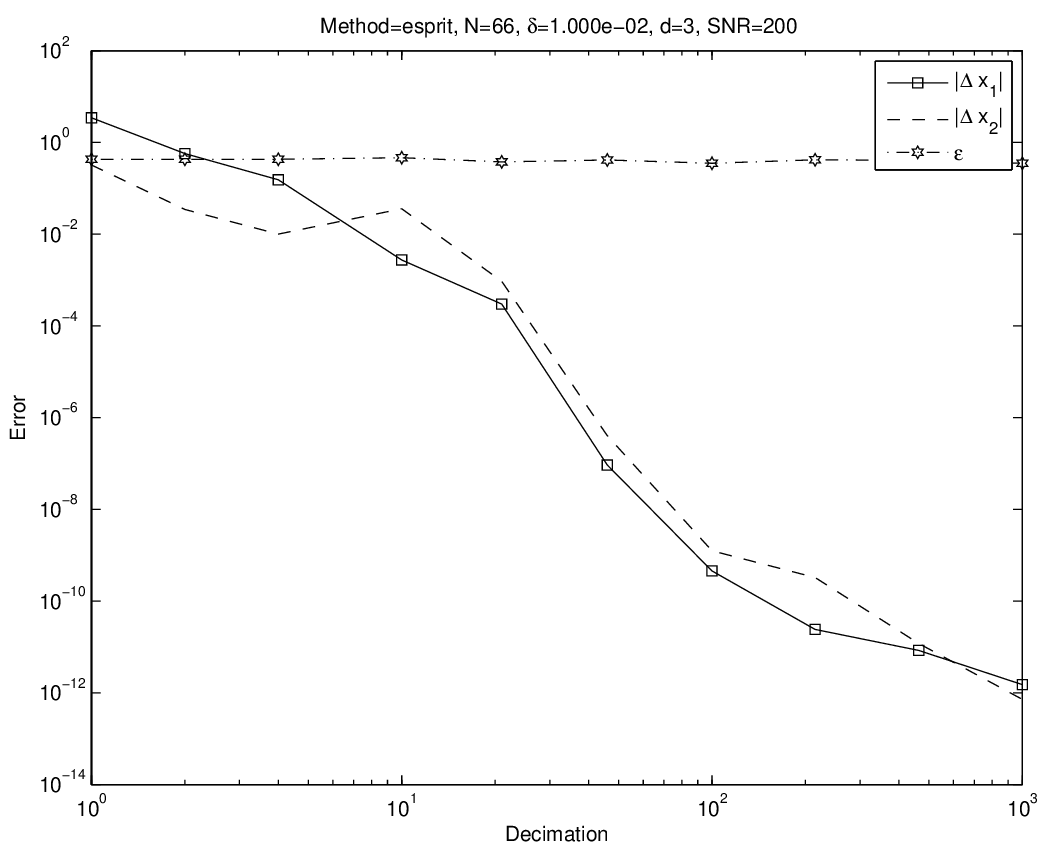}

}\subfloat[LS, $\ell_{j}=3$]{\includegraphics[width=0.45\columnwidth]{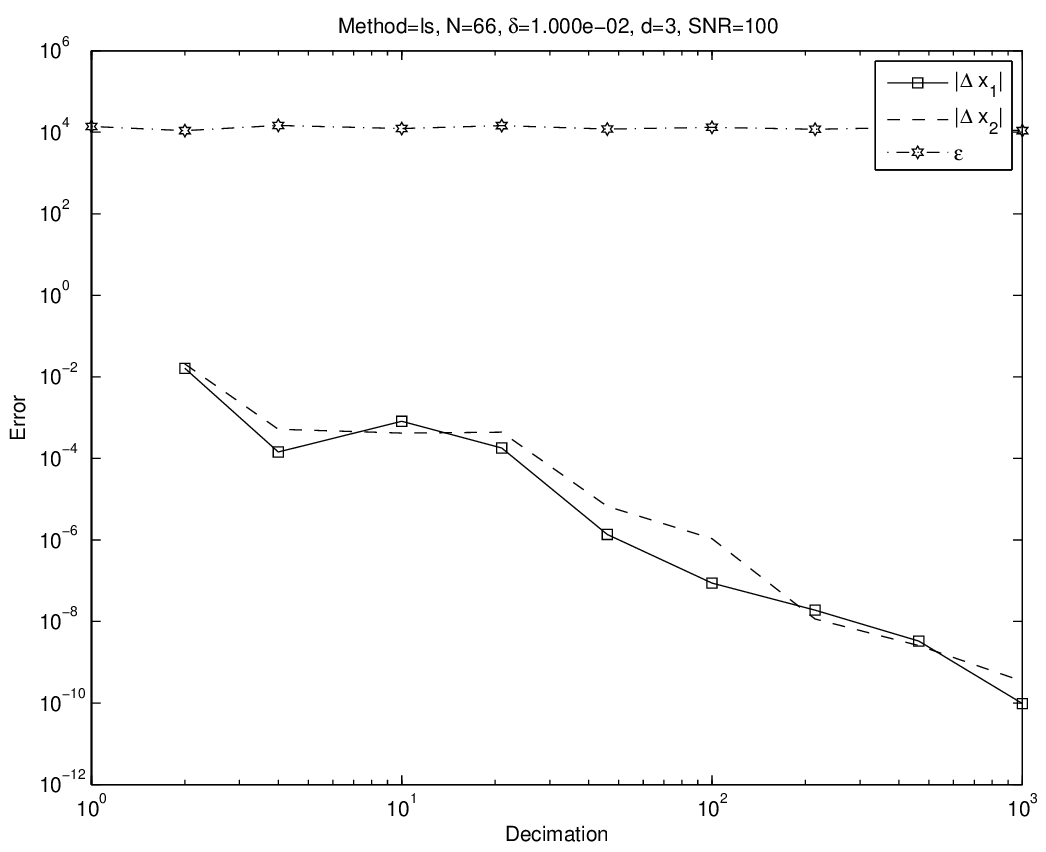}}

\caption[Decimation, fixed number of measurements]{Reconstruction error as a function of the decimation with fixed number
of measurements ($\protect\ns=66$). The signal has two nodes with
distance $\minsep=10^{-2}$ between each other. Notice that ESPRIT
requires significantly higher Signal-to-Noise Ratio in order to achieve
the same performance as LS.}
\label{fig:fixed-number-meas}
\end{figure}

In the first experiment, we fixed the number of measurements to be
66, and changed the decimation parameter $\df$, while keeping the
noise level constant. The accuracy of recovery increased with $\df$
\textendash{} see \prettyref{fig:fixed-number-meas}.

\begin{figure}
\subfloat[ESPRIT, $\ell_{j}=2$]{\includegraphics[width=0.45\columnwidth]{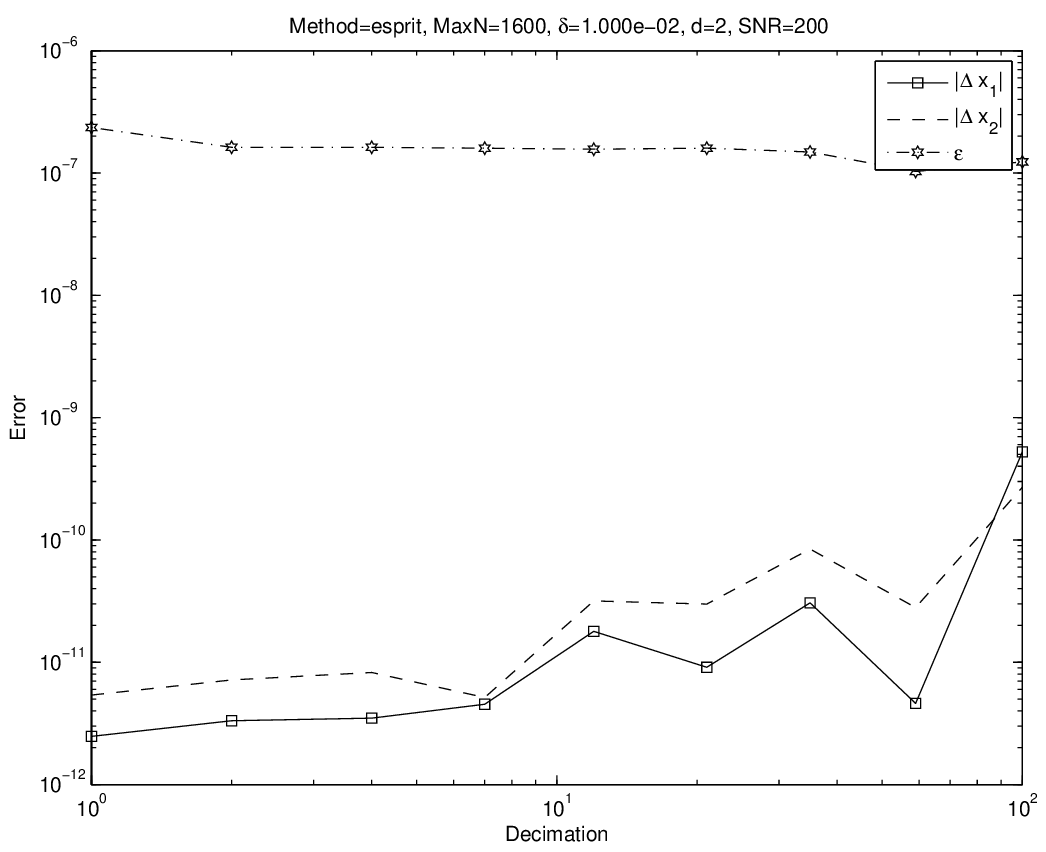}

} \subfloat[LS, $\ell_{j}=2$]{\includegraphics[width=0.45\columnwidth]{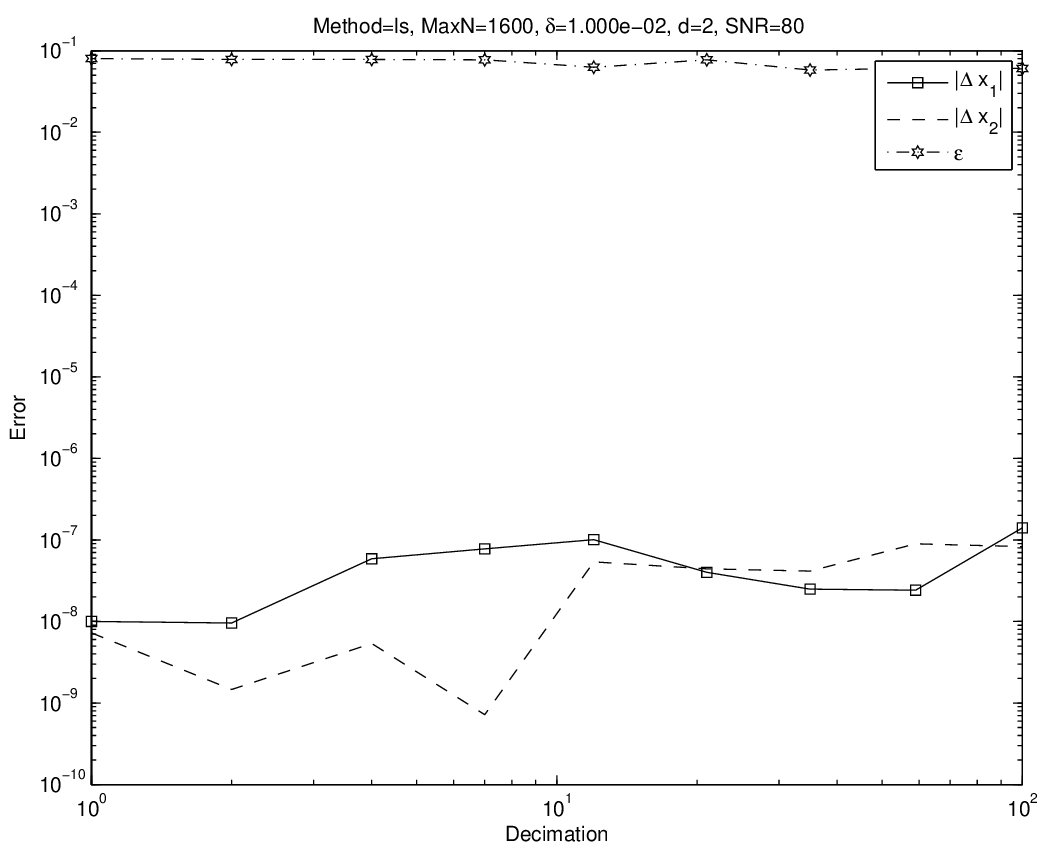}}

\caption[Decimation, reducing number of measurements]{Reconstruction error as a function of the decimation, reducing number
of measurements from $\protect\ns=1600$ to $\protect\ns=16$. The
signal has two nodes with distance $\minsep=10^{-2}$ between each
other. The reconstruction accuracy remains almost constant.}
\label{fig:full-decimation}
\end{figure}

In the second experiment, we fixed the highest available measurement
to be $\ns=1600$, and changed the decimation from $\df=1$ to $\df=100$
(thereby reducing the number of measurements from $1600$ to just
$16$). The accuracy of recovery stayed relatively constant \textendash{}
see \prettyref{fig:full-decimation}. Such a reduction leads to a
corresponding decrease in the running time, since for instance the
SVD computation in ESPRIT takes $O\left(\ns^{2}\nparams\right)$.

\section{\label{sec:Relation-to-existing}Relation to existing work}

Majority of the existing works in the literature consider the first
order Prony system \eqref{eq:basic-prony}. Specializing the results
of the present paper to this special case, we have the following result.
\begin{thm}
\label{thm:specialization-to-simple-prony}Consider the system \eqref{eq:polynomial-prony}
with $\ell_{1}=\dots=\ell_{\np}=1$, and with a-priori bounds as elaborated
in \prettyref{sec:Conditioning}. 

\begin{enumerate}
\item For $\ns\minsep\gg1$ and for $j=1,\dots,\np$ we have (here $A=\sum_{m=1}^{\np}\left|\jc_{j}\right|$)
\begin{eqnarray*}
\cn_{2j-1,\ns} & \lessapprox & A,\\
\cn_{2j,\ns} & \lessapprox & \frac{A}{\left|\jc_{j}\right|}\cdot\frac{1}{\ns}.
\end{eqnarray*}
\item If, on the other hand, all the nodes form a cluster, i.e. $\ns\minsep^{*}<2\pi\np$,
then
\begin{eqnarray*}
\cn_{2j-1}^{\left(\df^{*}\right)} & \lessapprox & \frac{1}{\left(\ns\minsep^{\left(j\right)}\right)^{2\np}},\\
\cn_{2j}^{\left(\df^{*}\right)} & \lessapprox & \frac{1}{\left|\jc_{j}\right|}\cdot\frac{\minsep^{\left(j\right)}}{\left(\ns\minsep^{\left(j\right)}\right)^{2\np}}.
\end{eqnarray*}
\end{enumerate}
\end{thm}
In his influential paper \cite{donoho_superresolution_1992}, Donoho
gave bounds for noise amplification (modulus of continuity $\Lambda$)
for recovery of signed measures from their continuous spectra of width
$\Omega$ on a lattice with step size $\Delta$ in the superresolution
setting $\Omega\Delta\ll\pi$. The ratio $\frac{1}{\Delta\Omega}$ is called the ``super-resolution
factor'' (SRF). If the measure has at most $\ell$
nonzero coefficients\footnote{The original paper considers the ``sparse clumps'' model, where $\ell$ is understood as the density of spikes per unit interval. For our purposes it is sufficient to consider just the ``sparse'' model.}, then $\Lambda$
is shown to increase at least as $\approx\left(\frac{1}{\Delta}\right)^{2\ell-1}$
and at most as $\approx\left(\frac{1}{\Delta}\right)^{2\ell+1}$.
When $\Delta\to0$, the lower bound effectively scales as $\left(SRF\right)^{2\ell-1}$, and the same scaling was recently shown to hold also for the upper bound by Demanet and Nguyen \cite{demanet_recoverability_2014}.

No practical way to achieve the above bounds have been proposed, however,
 recent works of Candes and Fernandez-Granda \cite{candes_super-resolution_2013,candes_towards_2014,fernandez-granda_support_2013}
showed that under an additional assumption of node separation (effectively
putting $\ell=1$ above) a stable recovery via total variation (TV) minimization
is possible, both for the $\ell_{1}$-norm and for the locations of
the spikes. Additional recent works \cite{duval_exact_2014,tang_near_2015} explore penalized TV approaches and provide similar stability estimates under various assumptions.    

To express the above setting in the notations of this paper, we identify
$\Delta$ with $\minsep$, $\Omega$ with $\ns$ and $\ell$ with $\np$,
and put $\ell_{j}=1$. After this identification, the second part
of \prettyref{thm:specialization-to-simple-prony} gives an upper
bound for the modulus of continuity of the order $\left(SRF\right)^{2\ell}$,
which is slightly worse than the estimates in \cite{demanet_recoverability_2014,donoho_superresolution_1992}.
Our setting is more general however, as the spikes are not assumed to lie on a grid. Furthermore, we also provide perturbation bounds for the locations of
the spikes in terms of the super-resolution factor. 

In a recent paper \cite{moitra_threshold_2014} the authors observed
a phase transition for the (unstructured) condition number of Vandermonde
matrices, a clear analogy with our results (note that in addition
to a similar phase transition, our estimates also predict an exponential
increase w.r.t $\nparams$ in the condition number, see \prettyref{subsec:proof-undecimated-refined}). In another related work, Demanet and Townsend \cite{demanet_stable_2016} studied the problem of polynomial extrapolation of  analytic functions, and they showed two different stability regimes, depending on the number of samples of the function -- similar to what we have described in this paper. It would be highly interesting to relate these results to each other.

A method very similar to decimation, called ``subspace shifting'',
or interleaving, was proposed by Maravic \& Vetterli in \cite{maravic2005sar}
in the context of analyzing performance of Finite Rate of Innovation
(FRI) sampling in the presence of noise. Their idea was to interleave
the rows of the Hankel matrix used in subspace estimation methods,
effectively increasing the separation of closely spaced nodes. They
confirmed this idea with numerical experiments. The results of our
paper can be considered as a theoretical justification of their approach,
and its extension to the more general system \eqref{eq:polynomial-prony}.

In statistical signal estimation, the Cramer-Rao Lower Bound (CRB)
gives a lower bound for the variance of any unbiased estimator, see
\cite{kay1993fundamentals}). In \cite{kusuma_accuracy_2009} the
authors only prove the CRB estimates for $\np=1,2$ and $\ns\gg1$,
for the system \eqref{eq:basic-prony}. On the other hand, the authors
of \cite{badeau_cramerrao_2008} consider the more general system
\eqref{eq:polynomial-prony} (called PACE model), and derive asymptotic
estimates for $\ns\gg1$. These results are qualitatively similar
to our \prettyref{thm:prony-square-cond} and \prettyref{thm:conditioning-full-stable}.
Obviously our results are different in nature from the CRB, but nevertheless
the stated similarity is worth investigating further. Generalized
ESPRIT is shown to asymptotically attain the CRB for $\ns\to\infty$.

The effect of oversampling for FRI signals was also studied in \cite{ben-haim_performance_2012},
where they showed that it can improve performance by several orders
of magnitude - a conclusion which is certainly consistent with our
\prettyref{thm:conditioning-full-stable}.

Stability analysis of Approximate Prony method, carried out by the
authors of \cite{peter_nonlinear_2011,potts_parameter_2010}, suggests
an increase in recovery error for the linear coefficients $\jc_{j}$,
again consistent with our results (see \cite{batenkov_accuracy_2013}
for further details).

Performance analysis of MUSIC in another recent paper
\cite{liao_music_2014} (see also a recent preprint \cite{fannjiang_compressive_2016} regarding ESPRIT) suggests that it can resolve arbitrarily close
frequencies below $\ns^{-1}$ for sufficiently small noise - compare
this with \prettyref{thm:conditioning-full-stable}, which shows that
the sensitivity indeed does not depend on the node separation.

The method of Filbir et. al \cite{filbir_problem_2012} solves the
system \eqref{eq:basic-prony} via constructing a certain orthogonal
polynomial on the unit circle. Their perturbation analysis gives an
error in the nodes of the order of $\sqrt{\frac{\log\ns}{\ns}}$. Also, localized kernel methods were recently shown to provide stable  estimation of instantaneous frequencies, under minimal separation assumption \cite{chui_signal_2015}. 

Decimation has recently appeared in zooming methods such as ZMUSIC
\cite{kia_high-resolution_2007} and zoom-ESPRIT \cite{kim_high-resolution_2013}
for reducing computational complexity and memory requirements for
estimating frequencies in a specified range. Experiments show also
improvement in accuracy of the zooming techniques w.r.t to their regular
counterparts, thus it would be interesting to see whether an analysis
similar to ours can be applied also in these cases.

A variant of decimation for Prony systems, called ``arithmetic progression
sampling'' (APS) and described in detail in \cite{sidi2003practical},
was shown by Sidi to enhance substantially both the convergence acceleration
and numerical stability properties of generalizations of the Richardson
extrapolation process. It would be interesting to make this connection
more elaborate and precise.

A kind of ``stochastic decimation'' (randomized arithmetic progression
sampling) was recently used by Kaltofen et.al for outlier removal
in sparse model synthesis and interpolation \cite{kaltofen_cleaning-up_2014}.

\section{\label{sec:Some-future-directions}Some future directions}

This paper is a part of a continuing research effort, investigating
the applicability of algebraic methods to signal reconstruction problems
\cite{akinshin_accuracy_2015,batenkov_prony_2014,batenkov_complete,batenkov_accurate,batenkov_accuracy_2013, batenkov_geometry_2014,batenkov_local_2013,batenkov_algebraic_2012,ettinger_linear_2008,sarig_signal_2008}.
Some of its findings were initially reported in \cite{batenkov_algebraic_2013}.
Building upon the presented ideas, we have recently proposed a novel
``decimated homotopy'' algorithm, which has been shown to achieve the accuracy specified in \prettyref{cor:super-resolution}, and outperform
state of the art methods such as ESPRIT in the near-colliding setting
\cite{batenkov_prony_2014,batenkov_accurate}. Another extension
of this work is reported in \cite{akinshin_accuracy_2015}, providing
tight global bounds (opposed to the first-order situation of this
paper) for the accuracy of cluster recovery. Decimation also played
a major role in our recent proposed algorithm for resolving the Gibbs
phenomenon \cite{batenkov_complete}. 

The numerical analysis of Prony systems in an important topic for
further investigations. For instance, the bounds of \prettyref{thm:conditioning-full-stable}
are valid for the noise model \eqref{eq:noise-model}. However, in
some applications such as \cite{batenkov_complete}, a more appropriate
assumption is
\[
\frac{\left|\Delta\meas\right|}{\left|\meas\right|}\leqslant\rho k^{-1},
\]
for some fixed $\rho$. In general, ``semi-global'' analysis is
required in this and similar settings, and we leave this for a future
publication (cf. \cite{akinshin_accuracy_2015}).

An important open question connected with stable solution of Prony
systems is how to detect the near-singular situations, and choose
the problem structure vector $\str=\left(\ell_{1},\dots,\ell_{\np}\right)$
in an optimal way. One possible approach might involve symbolic-numeric
techniques for polynomial systems, combined with analysis of the singularities
of the mapping $\fwm$ (\cite{batenkov_geometry_2014,batenkov_local_2013}). 

Under our assumption of a single cluster, decimation appears to provide
near-optimal conditioning with respect to the number of samples $\ns$.
While theoretical justification of this optimality would be desirable,
a more important goal is to provide optimal solution when only \emph{some}
of the nodes form a cluster.

\section{Acknowledgements}
We would like to thank H.~Mhaskar for useful suggestions regarding the manuscript.

\section*{References}

\bibliographystyle{abbrvnat}
\bibliography{decimated-prony-new}

\end{document}